\documentclass[a4paper,12pt]{amsart}
\usepackage{hyperref,fancyhdr,mathrsfs,amsmath,amscd,amsthm,amsfonts,latexsym,amssymb}

\DeclareMathOperator{\trace}{trace}

\DeclareMathOperator{\dif}{d}

\newcommand{\Cal}{\mathcal{C}}
\renewcommand{\H}{\mathscr{H}}
\newcommand{\V}{\mathscr{V}}
\newcommand{\F}{\mathscr{F}}
\newcommand{\Fa}{\mathcal{F}}
\newcommand{\J}{\mathcal{J}}

\def \a{\alpha}

\def \G{\Gamma}
\def \l{\lambda}

\def \phi{\varphi}
\def \Phi{\varPhi}
\def \p{\pi}

\def \t{\tau}

\def \R{\mathbb{R}}

\def \C{\mathbb{C}\,}

\def\widecheckg{g^{\hspace*{-2.5pt}\vbox to 5pt{\hbox to
0pt{\LARGE$\check{}$}}}\hspace*{2pt}}

\def\widecheckl{\lambda^{\hspace*{-3.5pt}\vbox to 8pt{\hbox to
0pt{\LARGE$\check{}$}}}\hspace*{2pt}}

\begin{document}

\title{Harmonic morphisms between Weyl spaces and twistorial maps II}
\author{Eric~Loubeau and Radu~Pantilie\;\dag}
\thanks{\dag\;Gratefully acknowledges that this work was partially supported by a
C.N.C.S.I.S. grant no.\ GR202/19.09.2006.}
\email{\href{mailto:Eric.Loubeau@univ-brest.fr}{Eric.Loubeau@univ-brest.fr},
       \href{mailto:Radu.Pantilie@imar.ro}{Radu.Pantilie@imar.ro}}
\address{E.~Loubeau, D\'epartement de Math\'ematiques, Laboratoire C.N.R.S. U.M.R. 6205,
Universit\'e de Bretagne Occidentale, 6, Avenue Victor Le Gorgeu, CS 93837,
29238 Brest Cedex 3, France}
\address{R.~Pantilie, Institutul de Matematic\u a ``Simion Stoilow'' al Academiei Rom\^ane,
C.P. 1-764, 014700, Bucure\c sti, Rom\^ania}

\newtheorem{thm}{Theorem}[section]
\newtheorem{lem}[thm]{Lemma}
\newtheorem{cor}[thm]{Corollary}
\newtheorem{prop}[thm]{Proposition}

\theoremstyle{definition}

\newtheorem{defn}[thm]{Definition}
\newtheorem{rem}[thm]{Remark}
\newtheorem{exm}[thm]{Example}

\numberwithin{equation}{section}

\maketitle
\thispagestyle{empty}
\vspace{-3mm}
\begin{center}
\emph{This paper is dedicated to the memory of James Eells.}
\end{center}
\vspace{1mm}

\section*{Abstract}
\begin{quote}
{\footnotesize  We define, on smooth manifolds, the notions of almost twistorial structure
and twistorial map, thus providing a unified framework for all known examples of twistor spaces.
The condition of being harmonic morphisms naturally appears among the geometric properties of
submersive twistorial maps between low-dimensional Weyl spaces endowed with a nonintegrable
almost twistorial structure due to Eells and Salamon. This leads to the twistorial characterisation
of harmonic morphisms between Weyl spaces of dimensions four and three. Also, we give a thorough
description of the twistorial maps with one-dimensional fibres from four-dimensional Weyl spaces
endowed with the almost twistorial structure of Eells and Salamon.}
\end{quote}

\section*{Introduction}

\indent
A predominant theme in the theory of harmonic maps is their relationship to
holomorphicity.
Right from their inception, harmonic maps have been recognised by Eells
and Sampson to include holomorphic maps between K\"ahler
manifolds, and a few years later, Lichnerowicz distinguished two
sub-classes of almost Hermitian manifolds for which this also holds.\\
\indent
A far harder task is to find conditions forcing harmonic maps to be
holomorphic. The clearest situation is between two-spheres
since harmonicity then implies conformality, but, in general, additional
conditions are required, on the map or on curvature.\\
\indent
From a Riemann surface, or a complex projective space, into a compact irreducible Hermitian
symmetric space harmonic stability ensures holomorphicity \cite{BBdBR}\,, \cite{OhUd}\,.\\
\indent
In the broader context of compact K\"ahler manifolds, Siu \cite{Siu} used a
$\partial\overline{\partial}$-Bochner argument to prove, under strong negativity
of the target curvature, the holomorphicity of rank four harmonic maps,
from which he deduced the biholomorphicity of compact K\"ahler manifolds
of the same homotopy type.\\
\indent
Whilst these results are of great importance in understanding harmonic
maps, this scheme bears its own natural limits, if only because it
requires the presence of complex structures, with their topological
consequences.\\
\indent
To overcome this hurdle, one can replace the codomain with an adequate
bundle admitting a natural complex structure, such that harmonicity of
any map is given by holomorphicity of its lift to this bundle.\\
\indent
This strategy, which could be traced to Calabi and even Weierstrass,
was put into effect by Eells and Salamon \cite{EelSal}\,, who modified
a well-known twistor construction of Atiyah, Hitchin and Singer
(see Example \ref{exm:smoothtwiststr4}\,, below)
to define, in our terminology, an \emph{almost twistorial structure}
on any oriented four-dimensional Riemannian manifold (see Example \ref{exm:smoothtwiststr4'}\,).
This yields a bijective correspondence between conformal harmonic maps and holomorphic curves
(see Proposition \ref{prop:smoothtwistmap2to3'4'}\,).
These ideas were extensively pursued by Bryant \cite{Bry-twistor}\,,
and Burstall and Rawnsley \cite{BurRaw-LNM} for (even-dimensional) Riemannian symmetric spaces.\\
\indent
The success of this strategy leads naturally to considering lifts on
both the domain and the codomain, hence removing any need of pre-existing
almost complex structures. The objective is two-fold: firstly, show that the
existence of a holomorphic lift implies harmonicity and, secondly, find natural
conditions under which harmonic maps admit holomorphic lifts.\\
\indent
As holomorphic maps are closed under composition, it seems that harmonic morphisms
will have an important role in this programme.\\
\indent
Besides, the very conformal nature of their characterization makes Weyl
geometry an ideal framework for their study, but it also turns out to
provide examples of twistor spaces in dimension two, three and four.\\
\indent
In the complex category, a twistor, that is, a point of the twistor space of a (complex) manifold $M$,
determines a submanifold, possibly, with self-intersections, of $M$. For example, the twistor space of
an anti-self-dual complex-conformal four-dimensional manifold $(M^4,c)$ is the space of self-dual surfaces
of $(M^4,c)$\,; also, the twistor space of a three-dimensional Einstein--Weyl space $(M^3,c,D)$ is the space of
coisotropic surfaces of $(M^3,c)$ which are totally-geodesic with respect to $D$ (see \cite{PanWoo-sd}\,).
Similar considerations apply to the flat complex-conformal even-dimensional manifolds
and to the complex-Riemannian manifolds of constant curvature (see \cite{Pan-sp}\,).
Also, the space of isotropic geodesics of a complex-conformal manifold is, in a natural way,
a twistor space \cite{LeB-nullgeod}\,.\\
\indent
A complex analytic map $\phi:M\to N$ between manifolds endowed with twistorial structures is
twistorial if it maps (some of the) twistors on $M$ to twistors on $N$ (see \cite{PanWoo-sd}\,). \\
\indent
In this paper we extend the notions of almost twistorial structure and twistorial map, to the smooth category.
We show that, \emph{in the smooth category, a twistor on a manifold $M$ (endowed with a twistorial structure)
is a pair $(N,J)$ where $N$ is a submanifold of $M$, possibly, with self-intersections, and $J$ is a linear
CR-structure on the normal bundle of $N$ in $M$} (Remark \ref{rem:smoothtwiststr}\,).
These submanifolds may well just be points. For example,
the twistor space of a four-dimensional anti-self-dual conformal  manifold $(M^4,c)$ is formed of pairs
$(x,J_x)$ where $x\in M$ and $J_x$ is a positive orthogonal complex structure on $(T_xM,c_x)$\,. Similarly,
the twistor space $Z$ of a three-dimensional conformal manifold $(M^3,c)$ is a five-dimensional CR-manifold
consisting of the orthogonal nontrivial CR-structures on $(M^3,c)$ \cite{LeB-CR} (see Example \ref{exm:nullgeodesics}
and note that, if we assume $(M^3,c)$ real-analytic, then $Z$ is a real hypersurface,
endowed with the induced CR-structure, in the space of isotropic geodesics of the germ-unique
complexification of $(M^3,c)$\,).\\
\indent
On the other hand, the twistor space of a three-dimensional Einstein--Weyl space
(see Example \ref{exm:smoothtwiststr3} and Theorem \ref{thm:smoothtwist3}) is formed of pairs
$(\gamma,J)$ where $\gamma$ is a geodesic and $J$ is an orthogonal complex structure on the normal bundle
of $\gamma$\,. The natural generalizations of this example, to Riemannian manifolds of constant curvature,
are given in \cite{Pan-sp}\,.\\
\indent
As in the complex category, a (smooth) map $\phi:M\to N$ between manifolds endowed with
twistorial structures is twistorial if it maps twistors on $M$ to twistors on $N$
(Definition \ref{defn:smoothtwistmap}\,). It follows that twistorial maps naturally generalize
holomorphic maps. Examples of twistorial maps have been previously used to obtain
constructions of Einstein and anti-self-dual manifolds (see \cite{Cal-sds}\,, \cite{PanWoo-exm}\,, \cite{PanWoo-sd}
and the references therein).\\
\indent
In Section 1, we recall a few basic facts on harmonic morphisms between Weyl spaces. Section 2
is preparatory for Sections 3 and 4\,, where we give the definitions of almost twistorial structures
and twistorial maps, on smooth manifolds. In Sections 3 and 4\,, we also give
examples of almost twistorial structures and twistorial maps.\\
\indent
In the first part \cite{LouPan} of this work, we introduced the notion of harmonic morphism between (complex-)Weyl spaces
and we continued the study (initiated in \cite{PanWoo-sd}\,) of the relations between harmonic morphisms
and twistorial maps. As all of the known examples of complex analytic almost twistorial structures
have smooth versions, all of the main results of \cite{LouPan} have `real' versions.
But not all smooth almost twistorial structures come from complex analytic almost twistorial structures.
The first example of such an almost twistorial structure is the almost twistorial structure of Eells and Salamon,
mentioned above. We introduce similar almost twistorial structures
(Examples \ref{exm:smoothtwiststr3'} and \ref{exm:smoothtwiststr4'phi}\,)
with respect to which a map (between Weyl spaces of dimensions four and three) is twistorial
if and only if it is a harmonic morphism (Theorem \ref{thm:smoothtwistmap4'phi3'}(i)\,);
this implies that \emph{there exists a bijective correspondence between one-dimensional
foliations on $(M^4,c,D)$\,, which are locally defined by harmonic morphisms, and certain almost
CR-structures on the bundle of positive orthogonal complex structures on $(M^4,c)$}\,. Other
similar results are given in Section 5\,.\\
\indent
In Section 6\,, we also give the necessary and sufficient conditions for a map with one-dimensional
fibres from a four-dimensional Weyl space endowed with the almost twistorial structure of Eells and Salamon
to be twistorial (Theorem \ref{thm:smoothtwistmap4'to3''}\,). It follows that, all such twistorial maps
are harmonic morphisms. Also, a generalized monopole equation (Definition \ref{defn:pullbackmono}(1)\,)
is naturally involved.

\section{Harmonic morphisms}

\indent
In this section all the manifolds and maps are assumed to be smooth.\\
\indent
Let $M^m$ be a manifold of dimension $m$\,. If $m$ is even then
we denote by $L$ the line bundle associated to the frame bundle of $M^m$
through the morphism of Lie groups $\rho_m:{\rm GL}(m,\R)\to(0,\infty)$\,,
$\rho_m(a)=|\det a\,|^{1/m}$\,, $\bigl(a\in{\rm GL}(m,\R)\bigr)$\,;
obviously, $L$ is oriented. If $m$ is odd then we denote by $L$
the line bundle associated to the frame bundle of $M^m$ through
the morphism of Lie groups $\rho_m:{\rm GL}(m,\R)\to\R^*$\,,
$\rho_m(a)=(\det a)^{1/m}$\,, $\bigl(a\in{\rm GL}(m,\R)\bigr)$\,;
obviously, $L^*\otimes TM$ is an oriented vector bundle.
We say that $L$ is \emph{the line bundle of $M^m$} (cf.\ \cite{Cal-sds}\,).\\
\indent
Similar considerations apply to any vector bundle.\\
\indent
Let $\phi:M\to N$ be a submersion and let $\H$ be a distribution on $M$,
complementary to the fibres of $\phi$\,. Let $L_{\H}$ and $L_N$ be the line
bundles of $\H$ and $N$, respectively. As $L_{\H}^{\:2n}=\bigl(\Lambda^n(\H)\bigr)^2$
and $L_N^{\:2n}=\bigl(\Lambda^n(TN)\bigr)^2$, where $n=\dim N$, the differential
of $\phi$ induces a bundle map $\Lambda$ from $L_{\H}^{\:2}$ to $L_N^{\:2}$.
If $n$ is odd then we also have a bundle map $\l$ from $L_{\H}$ to $L_N$\,;
obviously, $\Lambda=\l^2$. Furthermore, if $n$ is odd, $\dif\!\phi$ and $\l$ induce
a bundle map from $L_{\H}^*\otimes\H$ to $L_N^*\otimes TN$, which will also be denoted
by $\dif\!\phi$\,; note that, $\dif\!\phi:L_{\H}^{\:*}\otimes\H\to L_N^{\:*}\otimes TN$ is
orientation preserving, on each fibre.\\
\indent
Let $c$ be a conformal structure on $M$; that is, $c$ is a section of $L^2\otimes\bigl(\odot^2T^*\!M\bigr)$
which is `positive-definite'; that is, for any positive section $s^2$ of $L^2$ we have
$c=s^2\otimes g_s$\,, where $g_s$ is a Riemannian metric on $M$;
then $g_s$ is a \emph{representative} of $c$\,. Therefore, $c$ corresponds to a Riemannian metric
on the vector bundle $L^*\otimes TM$ (see \cite{Cal-sds}\,).
Furthermore, $c$ corresponds to an injective vector bundle morphism $L^2\hookrightarrow\odot^2TM$
such that, at each $x\in M$, the positive bases of $L^2_x$ are mapped into the cone of positive
definite symmetric bilinear forms on $T_x^*M$. Obviously, $c$ also corresponds to a reduction of
the frame bundle of $M$ to ${\rm CO}(m,\R)$\,, where $m=\dim M$; the total space of the
reduction corresponding to $c$ is formed of the \emph{conformal frames} on $(M,c)$\,.\\
\indent
If $\dim M$ is odd then local sections of $L$ correspond
to \emph{oriented local representatives} of $c$\,; that is, (local) representatives of $c$\,,
on some oriented open set of $M$. Let $\H$ be a distribution on $M$. Then $c$
induces a conformal structure $c|_{\H}$ on $\H$ and, it follows that, we have
an isomorphism, which depends of $c$\,, between $L^{\,2}$ and $L_{\H}^{\:2}$\,.

\begin{defn}[cf.\ \cite{BaiWoo2}\,]
Let $(M,c_M)$ and $(N,c_N)$ be conformal manifolds. A map $\phi:(M,c_M)\to(N,c_N)$
is \emph{horizontally weakly conformal} if at each point $x\in M$, either $\dif\!\phi_x=0$ or
$\dif\!\phi_x|_{({\rm ker}\dif\!\phi_x)^{\perp}}$ is a conformal linear isomorphism from
$\bigl(({\rm ker}\dif\!\phi_x)^{\perp}, (c_M)_x|_{({\rm ker}\dif\!\phi_x)^{\perp}}\bigr)$
onto $\bigl(\,T_{\phi(x)}N,(c_N)_{\phi(x)}\bigr)$\,.
\end{defn}

\indent
Let $(M,c_M)$ and $(N,c_N)$ be conformal manifolds and let $\phi:(M,c_M)\to(N,c_N)$
be a horizontally conformal submersion; denote by $L_M$ and $L_N$ the line bundles
of $M$ and $N$, respectively. Let $\H=({\rm ker}\dif\!\phi)^{\perp}$ and denote by
$L_{\H}$ its line bundle. As $L_M^{\:2}=L_{\H}^{\:2}$\,, the bundle map $\Lambda$ defines
an isomorphism from $L_M^{\:2}$ to $\phi^*(L_N^{\:2})$ which is called the
\emph{square dilation of $\phi$}\,. We have
$\Lambda=\tfrac1n\,c(\dif\!\phi,\dif\!\phi)$\,, where $n=\dim N$ and $c$ is the
conformal structure on the bundle $T^*M\otimes\phi^*(TN)$ induced by $c_M$\,, $c_N$
and $\phi$\,. Therefore, if $\phi$ is horizontally weakly conformal then $\Lambda$
extends to a (smooth) section of ${\rm Hom}\bigl(L_M^{\:2},\phi^*(L_N^{\:2})\bigr)$
which is zero over the set of critical points of $\phi$\,.\\
\indent
If $\dim N$ is odd, the bundle map $\l$ from $L_{\H}$ to
$L_N$ is called the \emph{oriented dilation of $\phi$}\,. Note that,
$\dif\!\phi:(L_{\H}^{\:*}\otimes\H,c|_{\H})\to(L_N^{\:*}\otimes
TN,c_N)$ is an orientation preserving isometry, on each fibre.\\
\indent
Let $(M,c)$ be a conformal manifold. A connection $D$ on $M$ is \emph{confomal}
if $Dc=0$\,; equivalently, $D$ is the covariant derivation of a principal
connection on the bundle of conformal frames on $(M,c)$\,. If $D$ is
torsion-free then it is called a \emph{Weyl connection} and $(M,c,D)$
is a \emph{Weyl space}.

\begin{defn}[cf.\ \cite{BaiWoo2}\,]
(i) Let $(M,c,D)$ be a Weyl space.
A \emph{harmonic function}, on $(M,c,D)$\,, is a function $f$, (locally) defined on $M$,
such that $\trace_c(D\!\dif\!f)=0$\,.\\
\indent
(ii) A map $\phi:(M,c_M,D^M)\to(N,c_N,D^N)$ between Weyl spaces is a \emph{harmonic map}
if $\trace_c(D\!\dif\!\phi)=0$\,,
where $D$ is the connection on $\phi^*(TN)\otimes T^*M$ induced by $D^M$, $D^N$ and $\phi$\,.\\
\indent
(iii) A map $\phi:(M,c_M,D^M)\to(N,c_N,D^N)$ between Weyl spaces is a \emph{harmonic morphism}
if for any harmonic function $f:V\to\R$\,, on $(N,c_N,D^N)$\,, with $V$ an open set of $N$
such that $\phi^{-1}(V)$ is nonempty, $f\circ\phi:\phi^{-1}(V)\to\R$ is a harmonic function,
on $(M,c_M,D^M)$\,.
\end{defn}

\indent
Obviously, any harmonic function is a harmonic map and a harmonic morphism, if $\R$
is endowed with its conformal structure and canonical connection.\\
\indent
The following result is basic for the theory of harmonic morphisms (see \cite{LouPan}\,,
\cite{BaiWoo2}\,).

\begin{thm}
A map between Weyl spaces is a harmonic morphism if and only if it is a harmonic map
which is horizontally weakly conformal.
\end{thm}

\indent
We end this section by recalling the fundamental equation for horizontally conformal
submersions between Weyl spaces (\,\cite{LouPan}\,; cf.\ \cite{BaiWoo2}\,).

\begin{prop} \label{prop:fundamn}
Let $\phi:(M^m,c_M,D^M)\to(N^n,c_N,D^N)$ be a horizontally conformal submersion between
Weyl spaces; $m=\dim M$, $n=\dim N$.\\
\indent
Let $L$ be the line bundle of $M^m$ and let $\V={\rm ker}\dif\!\phi$\,, $\H=\V^{\perp}$\,; 
denote by $D$ the Weyl connection of $(M^m,c_M,\V)$\,.\\
\indent
Then we have the following equality of partial connections on $L^{m-n}$, over $\H$,
\begin{equation*}
(\H D)^{m-n}+\trace_{c_M}(D\!\dif\!\phi)^{\flat}=(\H D^M)^{m-2}\otimes(D^N)^{-(n-2)}\;,
\end{equation*}
where we have denoted by the same symbol $D^N$ and its pull-back by $\phi$\,.
\end{prop}

\section{Complex distributions}

\indent
Unless otherwise stated, all the manifolds and maps are assumed to smooth.

\begin{defn} \label{defn:complexdistr}
A \emph{complex distribution} on a manifold $M$ is a complex subbundle $\F$ of
$T^{\C}\!M$ such that $\dim(\F_x\cap\overline{\F}_x)$\,, $(x\in M)$\,, is constant.
If $\F\cap\overline{\F}$ is the zero bundle then $\F$ is an \emph{almost CR-structure}
on $M$, and $(M,\F)$ an \emph{almost CR-manifold}.
\end{defn}

\begin{exm} \label{exm:Fstrinduced}
1) Let $F$ be an almost $f$-structure on a manifold $M$; that is, $F$ is a
section of ${\rm End}(TM)$ such that $F^3+F=0$ \cite{Yan-fstr}\,. We denote by
$T^0M$\,,\,$T^{1,0}M$\,,\,$T^{0,1}M$ the eigendistributions of
$F^{\C}\in\G\bigl({\rm End}(T^{\C}\!M)\bigr)$ corresponding
to the eigenvalues $0$\,,\,${\rm i}$\,,\,$-{\rm i}$\,, respectively.
Then $T^0M\oplus T^{0,1}M$ is a complex distribution on $M$ and
$T^{0,1}M$ is an almost CR-structure on $M$. We say that $T^0M\oplus T^{0,1}M$ is
\emph{the complex distribution associated to $F$}; similarly, $T^{0,1}M$ is
\emph{the almost CR-structure associated to $F$}.\\
\indent
2) Let $M$ be endowed with a complex distribution $\F$ and let $N\subseteq M$
be a submanifold. Suppose that $\dim(T^{\C}_x\!N\cap\F_x)$ and
$\dim(T^{\C}_x\!N\cap\F_x\cap\overline{\F}_x)$\,,  $(x\in N)$\,, are constant.
Then $T^{\C}\!N\cap\F$ is a complex distribution on $N$ which we call
\emph{the complex distribution induced by $\F$ on $N$}.\\
\indent
In particular, if $N$ is a real hypersurface in the complex manifold $M$ then
$T^{\C}\!N\cap T^{0,1}M$ is an almost CR-structure on $N$.
\end{exm}

\indent
Next, we define the notion of holomorphic map between manifolds endowed with complex
distributions.

\begin{defn} \label{defn:Fholo}
Let $\F^M$ and $\F^N$ be complex distributions on $M$ and $N$, respectively.
A map $\phi:(M,\F^M)\to(N,\F^N)$ is \emph{holomorphic} if $\dif\!\phi(\F^M)\subseteq\F^N$.\\
\indent
A map $\phi:(M,F^M)\to(N,F^N)$ between manifolds endowed with $f$-structures
is \emph{holomorphic} if $\phi:(M,\F^M)\to(N,\F^N)$ is holomorphic where
$\F^M$ and $\F^N$ are the complex distributions, on $M$ and $N$, associated
to $F^M$ and $F^N$, respectively.
\end{defn}

\indent
Let $F$ be an almost $f$-structure on $M$ and $J$ an almost complex structure on $N$.
Then, a map $\phi:(M,F)\to(N,J)$ is holomorphic if and only if
$\dif\!\phi\circ F=J\circ\dif\!\phi$\,.

\begin{defn}
Let $M$ be a manifold. A subbundle $E$ of $T^{\C}\!M$ is \emph{integrable}
if for any sections $X$,\,$Y$ of $E$ we have that $[X,Y]$ is a section of $E$\,.\\
\indent
A \emph{CR-structure} is an integrable almost CR-structure; a \emph{CR-manifold} is a manifold endowed
with a CR-structure.\\
\indent
An almost $f$-structure is \emph{integrable} if its associated complex distribution
is integrable; an \emph{$f$-structure} is an integrable almost $f$-structure.
\end{defn}

\indent
Note that a complex distribution $\F^M$ on $M$ is integrable if and only if for any point
$x\in M$, there exists a holomorphic submersion $\phi:(U,\F^M|_U)\to(N,\F^N)$\,, from
some open neighbourhood $U$ of $x$, onto some CR-manifold $(N,\F^N)$\,, such that
${\rm ker}\dif\!\phi=\bigl(\F^M\cap\overline{\F^M}\,\bigr)|_U$\,.
If $U=M$ and $\phi$ has connected fibres then $\F^M$ is called \emph{simple}; then,
obviously, up to a CR-diffeomorphism, $\phi$ is unique;
we call it the \emph{holomorphic submersion corresponding to $\F^M$}.\\

\indent
The following definition will, also, be useful in the next section.

\begin{defn}
Let $(P,J)$ be a complex manifold and let $Q\subseteq P$ be a complex
submanifold. A \emph{holomorphic distribution (on $P$) along $Q$} is a holomorphic
subbundle of $(T^{1,0}P)|_Q$\,.\\
\indent
Let $(M,\F)$ be a manifold endowed with an integrable complex distribution
such that $TM=\F+\overline{\F}$. Let $N\subseteq M$ be a
submanifold on which $\F$ induces a complex distribution. A complex subbundle $E$ of
$\F|_N$ is called a \emph{holomorphic distribution (on $M$) along $N$} if, for any
$x\in N$, there exists an open set $U\subseteq M$, containing $x$\,, and a holomorphic
submersion $\phi$ from $(U,\F|_U)$ onto some complex manifold $(P,J)$\,, with
${\rm ker}\dif\!\phi=(\F\cap\overline{\F})|_U$\,, such that $\phi|_{U\cap N}$
is a holomorphic submersion from $U\cap N$ onto
a complex submanifold of $(P,J)$ and $\dif\!\phi|_{U\cap N}$ maps $E|_{U\cap N}$,
isomorphically on each fibre, onto a holomorphic distribution along $\phi(U\cap N)$
(that is, onto a holomorphic subbundle of $T^{1,0}P|_{\phi(U\cap N)}$\,);
if $N=M$ then we say that $E$ is a \emph{holomorphic distribution on $M$}.
\end{defn}

\indent
Let $F$ be a $f$-structure on a manifold $M$. Obviously, $T^{1,0}M$ is a holomorphic
distribution on $M$.\\

\indent
Let $E$ be a complex vector bundle over $M$ endowed with a complex linear connection $\nabla$
(that is, $\nabla$ is a complex-linear map $\G(E)\to\G\bigl({\rm Hom}_{\,\C}(T^{\C}\!M,E)\bigr)$
such that $\nabla(sf)=(\nabla s)f+s\underset{\C}{\otimes}\dif\!f$
for any section $s$ of $E$ and any complex-valued function $f$ on $M$; equivalently,
we have $\nabla^{\C}=\nabla\oplus\overline{\nabla}$\,,
with respect to the decomposition $E^{\,\C}=E\oplus\overline{E}$\,).\\
\indent
Let $\pi:{\rm Gr}_k(E)\to M$ be the bundle of complex vector subspaces of complex dimension $k$ of $E$;
we shall denote by $\H\,\bigl(\subseteq T({\rm Gr}_k(E))\bigr)$ the connection induced by $\nabla$ on ${\rm Gr}_k(E)$\,.
Note that, as ${\rm Gr}_k(E)$ is a bundle whose typical fibre is a complex manifold and its structural group
a complex Lie group acting holomorphically on the fibre,
we have $({\rm ker}\dif\!\pi)^{\C}=({\rm ker}\dif\!\pi)^{1,0}\oplus({\rm ker}\dif\!\pi)^{0,1}$\,.\\
\indent
The following result, which we do not imagine to be new, will be useful later on.
We omit the proof.

\begin{prop} \label{prop:autoparallel}
Let $p$ be a section of ${\rm Gr}_k(E)$ (equivalently, $p$ is a complex vector subbundle of complex
rank $k$ of $E$\,). Let $X\in T^{\C}_{x_0}\!M$ for some $x_0\in M$;
denote by $p_0=p(x_0)$\,.\\
\indent
{\rm (a)} The following assertions are equivalent:\\
\indent
\quad{\rm (i)} $(\dif\!p)^{\C}(X)\in\H^{\C}_{p_0}$\,.\\
\indent
\quad{\rm (ii)} $\nabla^{\C}_Xs\in p_0^{\C}$ for any section $s$ of $p^{\C}$.\\
\indent
{\rm (b)} The following assertions are equivalent:\\
\indent
\quad{\rm (iii)} $(\dif\!p)^{\C}(X)\in\H^{\C}_{p_0}\oplus({\rm ker}\dif\!\pi)_{p_0}^{0,1}$.\\
\indent
\quad{\rm (iv)} $\nabla_Xs\in p_0$ for any section $s$ of $p$.
\end{prop}

Note that, if $X$ is real (that is, $X\in TM$) then the assertions (i)\,,\,\ldots\,,\,(iv)\,,
of Proposition \ref{prop:autoparallel}\,, are equivalent.

\section{Almost twistorial structures}

\indent
In this section, we define, in the smooth category, the notion of (almost)
twistorial structure.

\begin{defn}[cf.\ \cite{Raw-ftwistor}\,,\,\cite{PanWoo-sd}\,] \label{defn:smoothtwiststr}
An \emph{almost twistorial structure}, on a manifold $M$, is a quadruple
$\t=(P,M,\p,\F)$ where $\p:P\to M$ is a locally trivial fibre space and $\F$ is a
complex distribution on $P$ which induces almost complex structures on each fibre of $\p$\,;
if $\F$ is induced by an almost $f$-structure $F$ then, also, $(P,M,\p,F)$
is called an almost twistorial structure.\\
\indent
The almost twistorial structure $\t=(P,M,\p,\F)$ is \emph{integrable} if $\F$ is integrable.
A \emph{twistorial structure} is an integrable almost twistorial structure; the leaf space
of $\F\cap\overline{\F}$ is called the \emph{twistor space} of $\t$.\\
\indent
A twistorial structure $\t=(P,M,\p,\F)$ is called \emph{simple} if $\F$ is simple; if $\t$ is simple
with $\phi:(P,\F)\to Z$ the corresponding holomorphic submersion then $\dif\!\phi(\F)$ is
a CR-structure on $Z$.
\end{defn}

\begin{rem} \label{rem:smoothtwiststr}
Let $\t=(P,M,\p,\F)$ be a twistorial structure and let $Z$ be its twistor space.
Each $z\in Z$ determines a pair $(N_z,J_z)$ where $N_z$ is a
submanifold of $M$, possibly, with self-intersections, and $J_z$ is a linear CR-structure
on the normal bundle of $N_z$.\\
\indent
Indeed, let $N_z=\phi^{-1}(z)$\,, where $\phi:P\to Z$ is the (continuous) projection whose fibres
are the leaves of $\F\cap\overline{\F}$.
Then $\p|_{N_z}:N_z\to M$ is an immersion. Also, let
$({\rm ker}\dif\!\p)^{0,1}=\F\cap({\rm ker}\dif\!\p)^{\C}$. Then the restriction to $N_z$
of the quotient of $\F$ through $TN_z\oplus({\rm ker}\dif\!\p)^{0,1}$ defines a linear CR-structure $J_z$ on the normal
bundle $(\p|_{N_z})^*(TM)/TN_z$ of $N_z$ in $M$.
\end{rem}

\indent
Now, we formulate the examples of almost twistorial structures with which we shall work.

\begin{exm} \label{exm:smoothtwiststr2}
If $\F$ is a complex distribution on $M$ then $(M,M,{\rm Id}_M,\F)$ is an almost
twistorial structure.\\
\indent
In particular, let $(M^2,c)$ be a two-dimensional oriented conformal manifold.
Then, as the identity component of ${\rm CO}(2,\R)$ is isomorphic to $\C^{\!*}$, there exists
a unique (almost) Hermitian structure $J$ on $(M^2,c)$ such that
if $X\in TM$ then $(X,JX)$ is a positively oriented frame on $M^2$; then
$(M,M,{\rm Id},J)$ is a twistorial structure. Similarly, any oriented vector bundle
of (real) rank two, endowed with a conformal structure, is a complex line bundle.
\end{exm}

\begin{exm}[\,\cite{LeB-nullgeod}\,] \label{exm:nullgeodesics}
Let $(M^m,c)$ be a conformal manifold of dimension $m\geq2$\,. Let $\p:P\to M$ be the
bundle of null direction on $M^m$. Obviously, $P$ is a bundle with typical fibre
$Q_{m-2}$ and structural group ${\rm PCO}(m,\C)$ where $Q_{m-2}$ is the
nondegenerate hyperquadric in $\C\!P^{m-1}$\,; in particular, $\dim P=3m-4$\,.\\
\indent
Endow $(M^m,c)$ with a Weyl connection $D$ and let $\H\subseteq TP$ be the
connection induced by $D$ on $P$. We denote by $\H^{1,0}$ the subbundle of
$\H^{\C}$ such that, at each $p\in P$,
the subspace $\H^{1,0}_p\subseteq\H_p^{\C}$ is the horizontal
lift of $p\subseteq T_{\p(p)}M$. Then $\F=({\rm ker}\dif\!\p)^{1,0}\oplus\H^{1,0}$
is a CR-structure on $P$. Moreover, $\F$ does not depend of $D$\,.\\
\indent
Then $\t=(P,M,\p,\F)$ is a twistorial structure on $M^m$. If $(M^m,c)$ is
real analytic then, locally, $P$ is a (real) submanifold 
of the space of null geodesics $\mathcal{N}$ of the complexification of $(M^m,c)$\,; moreover, 
the CR-structure $\F$ of $P$ is induced by the complex structure of $\mathcal{N}$\,.\\
\indent
Note that, if $m=2$ then $P$ is a $\mathbb{Z}_2$-covering space of $M^2$ which
is trivial if and only if $M^2$ is oriented. Assuming $M^2$ oriented we obtain that
$T^{1,0}P$ restricted to one of the components of $P$ defines the Hermitian structure
of Example \ref{exm:smoothtwiststr2}\,.
\end{exm}

\begin{exm}[\,\cite{Hit-complexmfds}\,,\,\cite{Raw-ftwistor}\,] \label{exm:smoothtwiststr3}
Let $(M^3,c,D)$ be a three-dimensional Weyl space. Let $\p:P\to M$ be the
bundle of nonzero skew-adjoint $f$-structures on $(M^3,c)$\,. Obviously, $P$ is
also the bundle of nonzero skew-adjoint $f$-structures on the oriented Riemannian
bundle $(L^*\otimes TM,c)$\,. Therefore, $P$ is isomorphic
to the sphere bundle of $(L^*\otimes TM,c)$\,. In particular, the typical fibre
and the structural group of $P$ are $\C\!P^1$ and ${\rm PGL}(2,\C)$\,, respectively.\\
\indent
We could also define $P$ as follows: firstly, note that, there exists a unique
oriented Riemannian structure, on the vector bundle $E$ of skew-adjoint
endomorphisms on $(M,c)$\,, with respect to which $[A,B]=A\times B$\,, for any $A,B\in E$\,;
then $P$ is the sphere bundle of $E$\,.\\
\indent
The bundle $P$ is also isomorphic with the bundle of oriented two-dimensional
subspaces on $M^3$. Therefore, there exists a bijective correspondence between
one-dimensional foliations on $M^3$, with oriented orthogonal complement, and
almost $f$-structures on $(M^3,c)$\,. Furthermore, under this
bijection, conformal one-dimensional foliations correspond to (integrable)
$f$-structures.\\
\indent
Let $k$ be a section of $L^*$. We define a conformal connection connection $\nabla$ on $(M,c)$
by $\nabla_XY=D_XY+\tfrac12\,kX\times Y$ for any vector fields $X$ and $Y$ on $M$.
We say that $\nabla$ is \emph{the connection associated to $D$ and $k$}.\\
\indent
Let $\H\subseteq TP$ be the connection induced by $\nabla$ on $P$. We denote by
$\H^0$, $\H^{1,0}$ the subbundles of $\H^{\C}$ such that, at each $p\in P$,
the subspaces $\H^0_p\,,\,\H^{1,0}_p\subseteq\H_p^{\C}$ are the horizontal
lifts of the eigenspaces of $p^{\C}\in{\rm End}(T_{\p(p)}^{\C}\!M)$ corresponding
to the eigenvalues $0$\,, ${\rm i}$\,, respectively.\\
\indent
We define the almost $f$-structure $\Fa$ on $P$ with respect to which
$T^0P=\H^0$ and $T^{1,0}P=({\rm ker}\dif\!\p)^{1,0}\oplus\H^{1,0}$.
Then $\t=(P,M,\p,\Fa)$ is an almost twistorial structure on $M$.
\end{exm}

\begin{rem}
Let $(M^3,c)$ be a three-dimensional conformal manifold. Let $\p:P\to M$ be the
bundle of nonzero skew-adjoint $f$-structures on $(M^3,c)$\,.\\
\indent
Let $D$ be a Weyl connection on $(M^3,c)$ and let $\nabla$ be a conformal connection on $(M^3,c)$\,.
The following assertions are equivalent:\\
\indent
\quad(i) $\nabla$ and $D$ induce, by applying the construction of Example \ref{exm:smoothtwiststr3}
(with $\H$ the connection on $P$ induced by $\nabla$ and $D$\,, respectively), the same almost $f$-structure on $P$.\\
\indent
\quad(ii) $\nabla$ and $D$ are projectively equivalent.\\
\indent
\quad(iii) There exists a section $k$ of $L^*$ such that $\nabla$ is the connection
associated to $D$ and $k$\,.
\end{rem}

\begin{exm}[cf.\ \cite{EelSal}\,] \label{exm:smoothtwiststr3'}
Under the same hypotheses as in Example \ref{exm:smoothtwiststr3}\,,
we define the almost $f$-structure $\Fa'$ on $P$ with respect to which we have
$T^0P=\H^0$ and $T^{1,0}P=({\rm ker}\dif\!\p)^{0,1}\oplus\H^{1,0}$.
Then $\t'=(P,M,\p,\Fa')$ is an almost twistorial structure on $M$.
\end{exm}

\begin{exm}[\,\cite{AtHiSi}\,] \label{exm:smoothtwiststr4}
Let $(M^4,c,D)$ be a four-dimensional oriented Weyl space. Let $\p:P\to M$ be the
bundle of positive orthogonal complex structures on $(M^4,c)$\,. Obviously, $P$ is
also the bundle of positive orthogonal complex structures on the oriented Riemannian
bundle $(L^*\otimes TM,c)$\,. Let $E$ be the adjoint bundle of $(L^*\otimes TM,c)$
and let $*_c$ be the involution of $E$ induced by the Hodge star-operator
of $(L^*\otimes M^4,c)$\,, under the isomorphism $E=\Lambda^2(L\otimes T^*M)$\,.
Then $E=E_+\oplus E_-$ where $E_{\pm}$ is the vector bundle, of rank three,
formed of the eigenvectors of $*_c$ corresponding to the eigenvalue $\pm1$\,.
There exists a unique oriented Riemannian structure $<\cdot,\cdot>$
on $E_{\pm}$ with respect to which $AB=-<A,B>{\rm Id}_{TM}\pm\,A\times B$
for any $A,B\in E_{\pm}$\,. It follows that $P$ is the sphere bundle of $E_+$\,.\\
\indent
Similarly to Example \ref{exm:smoothtwiststr3}\,, there exists a bijective correspondence
between two-dimensional distributions $\F$ on $M^4$, with oriented orthogonal complement,
and pairs $(J,K)$ of almost Hermitian structures on $(M^4,c)$\,, with $J$ positive and
$K$ negative, such that $J|_{\F^{\perp}}=K|_{\F^{\perp}}$\,.\\
\indent
Let $\H\subseteq TP$ be the connection induced by $D$ on $P$. We denote by $\H^{1,0}$
the subbundle of $\H^{\C}$ such that, at each $p\in P$,
the subspace $\H^{1,0}_p\subseteq\H_p^{\C}$ is the horizontal lift of the
eigenspace of $p^{\C}\in{\rm End}(T_{\p(p)}^{\C}\!M)$ corresponding to the
eigenvalue ${\rm i}$\,.
We define the almost complex structure $\J$ on $P$ with respect to which
$T^{1,0}P=({\rm ker}\dif\!\p)^{1,0}\oplus\H^{1,0}$\,.
Then $(P,M,\p,\J)$ is an almost twistorial structure on $M$.
\end{exm}

\begin{exm}[\,\cite{EelSal}\,] \label{exm:smoothtwiststr4'}
Let $(M^4,c,D)$ be a four-dimensional oriented Weyl space. With the same notations as
in Example \ref{exm:smoothtwiststr4}, let $\J'$ be the almost complex structure on $P$
with respect to which $T^{1,0}P=({\rm ker}\dif\!\p)^{0,1}\oplus\H^{1,0}$.
Then $\t'=(P,M,\p,\J')$ is an almost twistorial structure on $M$.
\end{exm}

\indent
Next, we recall the necessary and sufficient conditions for the integrability
of the almost twistorial structures of Examples \ref{exm:smoothtwiststr3} and
\ref{exm:smoothtwiststr4}\,.

\begin{thm}[\,\cite{Hit-complexmfds}\,,\,\cite{GauTod}\,] \label{thm:smoothtwist3}
Let $(M^3,c,D)$ be a three-dimensional Weyl space and let $k$ be a section of the dual
of the line bundle of $N^3$. Let $\t=(P,M,\p,\Fa)$
be the almost twistorial structure of Example \ref{exm:smoothtwiststr3}\,.
Let $\H\,(\subseteq TP)$ be the connection induced by $\nabla$ on $P$.\\
\indent
Then\/ $\t$ depends only of $(M^3,c,D)$ whilst $T^{1,0}P$ depends only of $(M^3,c)$\,.
Moreover, $\t$ is integrable if and only if $(M^3,c,D)$ is Einstein--Weyl.\\
\indent
Furthermore, if $(M^3,c,D)$ is Einstein--Weyl then we have the following.\\
\indent
\quad{\rm (i)} $\H^{1,0}$ restricted to any fibre $P_x=\C\!P^1$ of $P$ is a holomorphic
distribution, along $P_x$\,, isomorphic to $\mathcal{O}(2)$\,, $(x\in M)$\,.\\
\indent
\quad{\rm (ii)} Locally, the following four assertions are equivalent:\\
\indent
\quad\quad{\rm (ii1)} $\H^{1,0}$ is a holomorphic distribution on $P$.\\
\indent
\quad\quad{\rm (ii2)} $\H$ is integrable.\\
\indent
\quad\quad{\rm (ii3)} The connection induced by $\nabla$ on $L^*\otimes TM$ is flat.\\
\indent
\quad\quad{\rm (ii4)} The scalar curvature $s^D$ of $D$ satisfies $s^D=\tfrac32\,k^2$,
and $*_cDk=F^D$ where $F^D$ is the curvature form of the connection induced by $D$
on $L$\,.
\end{thm}

\begin{thm}[\,\cite{AtHiSi}\,, see \cite{Cal-F}\,] \label{thm:smoothtwist4}
Let $(M^4,c,D)$ be a four-dimensional oriented Weyl space and let $\t=(P,M,\p,\J)$
be the almost twistorial structure of Example \ref{exm:smoothtwiststr4}\,.
Let $\H\,(\subseteq TP)$ be the connection induced by $D$ on $P$.\\
\indent
Then\/ $\t$ depends only of $(M^4,c)$\,. Moreover,\/ $\t$ is integrable if and only
if $(M^4,c)$ is anti-self-dual.\\
\indent
Furthermore, if $(M^4,c)$ is anti-self-dual then we have the following.\\
\indent
\quad{\rm (i)} $\H^{1,0}$ restricted to any fibre $P_x=\C\!P^1$ of $P$ is a holomorphic
distribution, along $P_x$\,, isomorphic to $\mathcal{O}(1)\oplus\mathcal{O}(1)$\,,
$(x\in M)$\,.\\
\indent
\quad{\rm (ii)} $\H^{1,0}$ is a holomorphic distribution on $P$ if and only if,
locally, either $D$ is the Levi-Civita connection of an Einstein representative
of\/ $c$ or $D$ is the Obata connection of a hyper-Hermitian structure on $(M^4,c)$\,.
\end{thm}

\begin{rem}
1) Let $(M^3,c,D)$ be a three-dimensional Weyl space and let $k$ be a section of $L^*$,
where $L$ is the line bundle of $N^3$. Let $T^{1,0}P$ be defined as in Example \ref{exm:smoothtwiststr3}\,.\\
\indent
Then by associating to any $p\in P$ the eigenspace corresponding to $-{\rm i}$ we can identify
$P$ with the bundle of null directions on $(N^3,c)$\,. Furthermore, under this isomorphism, $T^{1,0}P$
corresponds to the CR-structure of Example \ref{exm:nullgeodesics}\,, associated to $(N^3,c)$\,.
Conversely, any conformal connection on $(N^3,c)$ having this property is the connection associated
to a Weyl connection on $(N^3,c)$ and a section of $L^*$.\\
\indent
2) Let $(M^3,c)$ be a three-dimensional conformal manifold and let $D$ be a Weyl connection
on $(M^3,c)$\,. Also, let $k$ be a section of the dual of the line bundle $L$ of $M^3$.
Then the almost twistorial structure $\t'$ of Example \ref{exm:smoothtwiststr3'}\,,
associated to $(M^3,c,D,k)$ is nonintegrable; that is, always not integrable
(cf.\ \cite{EelSal}\,; apply the integrability result proved in \cite{Pan-sp}
to show that $\t'$ is nonintegrable). Furthermore, if $D_1$ is another Weyl connection on $(M^3,c)$
and $k_1$ a section of $L^*$ such that the almost twistorial of
Example \ref{exm:smoothtwiststr3'}\,, associated to $(M^3,c,D_1,k_1)$\,,
is equal to $\t'$ then $D=D_1$ and $k=k_1$\,.\\
\indent
3) Let $(M^4,c)$ be a four-dimensional oriented conformal manifold and let $D$ be a Weyl connection
on $(M^4,c)$\,. Then the almost twistorial structure $\t'$ of Example \ref{exm:smoothtwiststr4'}\,,
associated to $(M^4,c,D)$ is nonintegrable (cf.\ \cite{EelSal}\,; apply the integrability result
proved in \cite{Pan-sp} to show that $\t'$ is nonintegrable).
Furthermore, if $D_1$ is another Weyl connection on $(M^4,c)$
such that the almost twistorial structure of Example \ref{exm:smoothtwiststr4'}\,,
associated to $(M^4,c,D_1)$\,, is equal to $\t'$ then $D=D_1$\,.
\end{rem}

\section{Twistorial maps} \label{twistmaps}

\indent
We start this section with the definition of twistorial maps (cf.\ \cite{Raw-ftwistor}\,,\,\cite{PanWoo-sd}\,).

\begin{defn} \label{defn:smoothtwistmap}
Let $\t_M=(P_M,M,\p_M,\F^M)$ and $\t_N=(P_N,N,\p_N,\F^N)$ be almost
twistorial structures and let $\phi:M\to N$ be a map. Suppose that
there exists a locally trivial fibre subspace $\p_{M,\phi}:P_{M,\phi}\to M$ of $\p_M:P_M\to M$
and a map $\Phi:P_{M,\phi}\to P_N$ with the properties:\\
\indent
\quad1)\; $\F^M$ induces a complex distribution $\F^{M,\phi}$ on $P_{M,\phi}$ and almost complex structures
on each fibre of $\p_{M,\phi}$ such that $\dif\!\p_M(\F^M_{p})=\dif\!\p_{M,\phi}(\F^{M,\phi}_{p})\,$,
for any $p\in P_{M,\phi}$\,.\\
\indent
\quad2)\; $\phi\circ\p_{M,\phi}=\p_N\circ\Phi$\,.\\
\indent
Then $\phi:(M,\t_M)\to(N,\t_N)$ is a \emph{twistorial map (with respect to $\Phi$)} if the map
$\Phi:(P_{M,\phi},\F^{M,\phi})\to(P_N,\F^N)$ is holomorphic. If $\F^{M,\phi}$ and $\F^N$ are simple complex distributions,
with $(P_{M,\phi},\F^{M,\phi})\to Z_{M,\phi}$ and $(P_N,\F^N)\to Z_N$, respectively, the corresponding holomorphic submersions onto
CR-manifolds, then $\Phi$ induces a holomorphic map $Z_{\phi}:Z_{M,\phi}\to Z_N$
which is called the \emph{twistorial representation} of $\phi$\,.
\end{defn}

\begin{rem}
With the same notations as in Definition \ref{defn:smoothtwiststr}\,, we have that
$\t_{M,\phi}=(P_{M,\phi},M,\p_{M,\phi},\F^{M,\phi})$ is an almost twistorial structure on $M$.
Obviously, $\t_{M,\phi}$ is integrable if $\t_M$ is integrable. Then
from (1) it follows that the twistor space $Z_{M,\phi}$ of $\t_{M,\phi}$ is a topological
subspace of the twistor space $Z_M$ of $\t_M$. Moreover, for any
$z\in Z_{M,\phi}$ the pair $(N_z,J_z)$ of Remark \ref{rem:smoothtwiststr} applied to $\t_{M,\phi}$
is equal to the pair determined by $z$ as a point of $Z_M$.\\
\indent
If $\t_M$ is simple then $\t_{M,\phi}$ is also simple and $Z_{M,\phi}$ is a submanifold of $Z_M$.
Moreover, if we denote by $\Cal_M$ and $\Cal_{M,\phi}$ the CR-structures of $Z_M$ and $Z_{M,\phi}$\,,
respectively, then $\Cal_{M,\phi}=\Cal_M\cap T^{\C}\!Z_{M,\phi}$\,.
\end{rem}

\indent
Next, we give examples of twistorial maps between manifolds endowed with the almost twistorial
structures of Examples \ref{exm:smoothtwiststr2}\,,\,\ref{exm:smoothtwiststr3} or
\ref{exm:smoothtwiststr4}\,.

\begin{exm} \label{exm:smoothtwistmap2to2}
Let $(M^2,c_M)$ and $(N^2,c_N)$ be two-dimensional oriented conformal manifolds.
Let $\t_M=(M,M,{\rm Id}_M,J^M)$ and $\t_N=(N,N,{\rm Id}_N,J^N)$ be the twistorial
structures of Example \ref{exm:smoothtwiststr2} associated to $(M^2,c_M)$ and
$(N^2,c_N)$\,, respectively.\\
\indent
Let $\phi:M^2\to N^2$ be a map. Obviously, $\phi:(M^2,\t_M)\to(N^2,\t_N)$ is
twistorial (with respect to $\phi$) if and only if $\phi:(M^2,J^M)\to(N^2,J^N)$
is holomorphic.
\end{exm}

\indent
The following three examples are, essentially, due to \cite{Hit-complexmfds}\,, \cite{Raw-ftwistor}.

\begin{exm} \label{exm:smoothtwistmap1to3}
Let $(M^3,c,D)$ be a three-dimensional Weyl space and let $N^1$ be a
(connected) one-dimensional submanifold of $M^3$ such that $TN^{\perp}$ is orientable.\\
\indent
Let $\t_M=(P,M,\p,\Fa)$ be the almost twistorial
structure of Example \ref{exm:smoothtwiststr3}\,, associated to $(M^3,c,D)$\,.
Also, endow $N^1$ with the trivial $f$-structure and let $\t_N=(N,N,{\rm Id}_N,0)$ be the
twistorial structure given by $(N^1,0)$\,, as in Example \ref{exm:smoothtwiststr2}\,.\\
\indent
Let $p_N$ be the section of $P^5$ over $N^1$ such that, $(p_N)_x|_{T_xN}=0$\,, for any $x\in N$.\\
\indent
The following assertions are equivalent:\\
\indent
\quad(i) $(N^1,\t_N)\hookrightarrow(M^3,\t_M)$ is twistorial (with respect to $p_N$).\\
\indent
\quad(ii) $N^1$ is a geodesic of $D$\,.
\end{exm}

\begin{exm} \label{exm:smoothtwistmap2to3}
Let $(M^3,c,D)$ be a three-dimensional Weyl space and let $N^2$ be an
oriented surface in $M^3$. Let $\t_M=(P,M,\p,\Fa)$ be the almost twistorial structure
of Example \ref{exm:smoothtwiststr3}\,, associated to $(M^3,c,D)$\,, and
let $\t_N=(N,N,{\rm Id}_N,J)$ be the twistorial structure of
Example \ref{exm:smoothtwiststr2}\,, associated to $(N^2,c|_N)$\,.\\
\indent
Let $p_N$ be the section of $P^5$ over $N^2$ such that, $(p_N)_x|_{T_xN}=J_x$\,, for any $x\in N$.\\
\indent
The following assertions are equivalent:\\
\indent
\quad(i) $(N^2,\t_N)\hookrightarrow(M^3,\t_M)$ is twistorial (with respect to $p_N$).\\
\indent
\quad(ii) $p_N:(N^2,J)\to(P^5,\Fa)$ is holomorphic.\\
\indent
\quad(iii) $N^2$ is a totally umbilical submanifold of $(M^3,c)$\,.\\
\indent
It follows that the equivalence (i)$\iff$(iii) also holds for a surface in a conformal manifold
endowed with the twistorial structure of Example \ref{exm:nullgeodesics}\,.
\end{exm}

\begin{exm} \label{exm:smoothtwistmap3to2}
Let $(M^3,c_M,D)$ be a three-dimensional Weyl space and let $(N^2,c_N)$ be a
two-dimensional oriented conformal manifold.
Let $\t_M=(P,M,\p,\Fa)$ be the almost twistorial structure
of Example \ref{exm:smoothtwiststr3}\,, associated to $(M^3,c_M,D)$\,,
and let $\t_N=(N,N,{\rm Id}_N,J)$ be the twistorial structure of
Example \ref{exm:smoothtwiststr2}\,, associated to $(N^2,c_N)$\,.\\
\indent
Let $\phi:M^3\to N^2$ be a submersion and let $F^{\phi}$ be the almost $f$-structure
on $(M^3,c)$ determined by ${\rm ker}\dif\!\phi$ and the orientation
of $N^2$ (see Example \ref{exm:smoothtwiststr3}\,).
Define $\Phi=\phi\circ{p_{\phi}}^{-1}:p_{\phi}(M)\to N$ where
$p_{\phi}$ is the section of $P$ corresponding to $F^{\phi}$. Note that, $p_{\phi}$ restricted
to any fibre $\phi^{-1}(y)$ is equal to $p_{\phi^{-1}(y)}$ of Example \ref{exm:smoothtwistmap2to3}\,,
$(y\in\phi(M)\,)$\,.\\
\indent
The following assertions are equivalent:\\
\indent
\quad(i) $\phi:(M^3,\t_M)\to(N^2,\t_N)$ is twistorial (with respect to $\Phi$).\\
\indent
\quad(ii) $p_{\phi}:(M^3,F^{\phi})\to(P^5,\Fa)$ is holomorphic and
$\phi:(M^3,c_M)\to(N^2,c_N)$ is horizontally conformal.\\
\indent
\quad(iii) The fibres of $\Phi$ are tangent to the connection induced by $D$ on $P^5$
and $\phi:(M^3,c_M)\to(N^2,c_N)$ is horizontally conformal.\\
\indent
\quad(iv) $\phi:(M^3,c_M,D)\to(N^2,c_N)$ is a horizontally conformal submersion
with geodesic fibres.
\end{exm}

\begin{exm}[\,\cite{EelSal}\,] \label{exm:smoothtwistmap2to4}
Let $(M^4,c)$ be a four-dimensional oriented conformal manifold and let $N^2$ be an
oriented surface in $M^4$. Let $\t_M=(P,M,\p,\J)$ be the almost twistorial structure
of Example \ref{exm:smoothtwiststr4}\,, associated to $(M^4,c)$\,, and
let $\t_N=(N,N,{\rm Id}_N,J)$ be the twistorial structure of
Example \ref{exm:smoothtwiststr2}\,, associated to $(N^2,c|_N)$\,.\\
\indent
Let $p_N$ be the section of $P^6$ over $N^2$ such that $(p_N)_x|_{T_xN}=J_x$\,, for any $x\in N$.\\
\indent
The following assertions are equivalent:\\
\indent
\quad(i) $(N^2,\t_N)\hookrightarrow(M^4,\t_M)$ is twistorial (with respect to $p_N$).\\
\indent
\quad(ii) $p_N:(N^2,J)\to(P^6,\J)$ is holomorphic.\\
\indent
If we endow $M^4$ with the opposite orientation then Example \ref{exm:smoothtwiststr4}
gives another almost twistorial structure $\widetilde{\t}_M=(\widetilde{P},M,\widetilde{\p},\widetilde{\J})$\,.\\
\indent
The following assertions are equivalent:\\
\indent
\quad(a) $(N^2,\t_N)\hookrightarrow(M^4,\t_M)$ and $(N^2,\t_N)\hookrightarrow(M^4,\widetilde{\t}_M)$
are twistorial.\\
\indent
\quad(b) $N^2$ is a totally umbilical submanifold of $(M^4,c)$\,.
\end{exm}

\begin{exm}[\,\cite{Woo-4d}\,] \label{exm:smoothtwistmap4to2}
Let $(M^4,c_M)$ and $(N^2,c_N)$ be oriented conformal manifolds of dimensions
four and two, respectively. Let $\t_M=(P,M,\p,\J)$ be the almost twistorial structure
of Example \ref{exm:smoothtwiststr4}\,, associated to $(M^4,c_M)$\,,
and let $\t_N=(N,N,{\rm Id}_N,J)$ be the twistorial structure of
Example \ref{exm:smoothtwiststr2}\,, associated to $(N^2,c_N)$\,.\\
\indent
Let $\phi:M^4\to N^2$ be a submersion and let $J^{\phi}$ be the positive almost Hermitian structure
on $(M^4,c_M)$ determined by ${\rm ker}\dif\!\phi$ and the orientation
of $N^2$ (see Example \ref{exm:smoothtwiststr4}\,).
Define $\Phi=\phi\circ{p_{\phi}}^{-1}:p_{\phi}(M)\to N$ where
$p_{\phi}$ is the section of $P$ corresponding to $J^{\phi}$. Note that, $p_{\phi}$ restricted
to any fibre $\phi^{-1}(y)$ is equal to $p_{\phi^{-1}(y)}$ of Example \ref{exm:smoothtwistmap2to4}\,,
$(y\in\phi(M)\,)$\,.\\
\indent
The following assertions are equivalent:\\
\indent
\quad(i) $\phi:(M^4,\t_M)\to(N^2,\t_N)$ is twistorial (with respect to $\Phi$).\\
\indent
\quad(ii) $p_{\phi}:(M^4,J^{\phi})\to(P^6,\J)$ is holomorphic and
$\phi:(M^4,c_M)\to(N^2,c_N)$ is horizontally conformal.\\
\indent
\quad(iii) $\phi:(M^4,c_M)\to(N^2,c_N)$ is horizontally conformal and its fibres are
twistorial, in the sense of Example \ref{exm:smoothtwistmap2to4}\,.\\
\indent
\quad(iv) $J^{\phi}$ is integrable and $\phi:(M^4,c_M)\to(N^2,c_N)$ is horizontally
conformal.
\end{exm}

\begin{exm}[\,\cite{Hit-complexmfds}\,,\,\cite{Cal-sds}\,] \label{exm:smoothtwistmap4to3}
Let $(M^4,c_M)$ be a four-dimensional oriented conformal manifold and let
$(N^3,c_N,D^N)$ be a three-dimensional Weyl space. Let $\t_M=(P_M,M,\p_M,\J)$ be
the almost twistorial structure of Example \ref{exm:smoothtwiststr4}\,, associated
to $(M^4,c_M)$\,, and let $\t_N=(P_N,N,\p_N,\Fa)$ be the almost twistorial structure
of Example \ref{exm:smoothtwiststr3}\,, associated to $(N^3,c_N,D^N)$\,.\\
\indent
Let $\phi:M^4\to N^3$ be a submersion. Let $\V={\rm ker}\dif\!\phi$ and $\H=\V^{\perp}$.
Then the orientation of $M^4$ corresponds to an isomorphism, which depends of $c_M$\,,
between $\V$ and the line bundle of $\H$. Therefore $(\V^*\otimes\H,c|_{\H})$ is
an oriented Riemannian vector bundle. We define $\Phi:P_M\to P_N$ by
\begin{equation*}
\Phi(p)=\tfrac{1}{||\dif\!\phi(V^*\otimes p(V))||}\,\dif\!\phi(V^*\otimes p(V))\;,
\end{equation*}
where $\{V\}$ is any basis of $\V_{\p_M(p)}$ and $\{V^*\}$ its dual basis, $(p\in P_M)$\,.\\
\indent
Let $I^{\H}$ be the $\V$-valued two-form on $\H$ defined by $I^{\H}(X,Y)=-\V[X,Y]$\,,
for any sections $X$ and $Y$ of $\H$. Then $*_{\H}I^{\H}$ is a horizontal one-form
on $M^4$, where $*_{\H}$ is the Hodge star-operator of $(\V^*\otimes\H,c|_{\H})$\,.
Denote by $D_+$ the Weyl connection on $(M^4,c_M)$ defined by
$D_+=D+*_{\H}I^{\H}$ where $D$ is the Weyl connection of $(M^4,c_M,\V)$ (see \cite{LouPan}\,).\\
\indent
The following assertions are equivalent:\\
\indent
\quad(i) $\phi:(M^4,\t_M)\to(N^3,\t_N)$ is twistorial (with respect to $\Phi$).\\
\indent
\quad(ii) $\phi:(M^4,c_M)\to(N^3,c_N)$ is horizontally-conformal and
$\phi^*(D^N)=\H D_+$ as partial connections on $\H$, over $\H$.
\end{exm}

\indent
Let $\phi:(M^4,c_M)\to N^3$ be a submersion from a four-dimensional oriented conformal manifold to a three-dimensional
manifold. Let $L$ be the line bundle of $N^3$. As the orientation of $M^4$ corresponds to an isomorphism
between $\V$ and $\phi^*(L)$\,, the pull-back by $\phi$ of any section of $L^*$ is a (vertical) one-form
on $M^4$. Similarly, if $E$ is a vector bundle over $N^3$, the pull-back by $\phi$ of any section of
$L^*\otimes E$ is a $\phi^*(E)$-valued one-form on $M^4$.\\
\indent
Next, we recall the following definitions (see \cite{PanWoo-sd} and the references therein).

\begin{defn} \label{defn:pullbackmono}
Let $P$ be a principal bundle on $N^3$ endowed with a (principal) connection $\G$ and let $A$ be a section
of $L^*\otimes{\rm Ad}P$\,.\\
\indent
1) Let $N^3$ be endowed with a conformal structure $c_N$ and a Weyl connection $D^N$.
The pair $(A,\G)$ is called a \emph{monopole on $(N^3,c_N,D^N)$} if
$$R=*_{\!_N}(D^N\otimes\nabla)(A)$$ where $R$ is the curvature form of $\G$ and $\nabla$ is
(the covariant derivative of) the connection induced by $\G$\/ on ${\rm Ad}P$.\\
\indent
2) Let $(M^4,c_M)$ be a four-dimensional oriented conformal manifold and let $\phi:M^4\to N^3$ be a submersion.
The connection $\widetilde{\G}$ on $\phi^*(P)$ defined by $$\widetilde{\G}=\phi^*(\G)+\phi^*(A)$$
is called \emph{the pull-back by $\phi$ of $(A,\G)$}\,.
\end{defn}

\indent
Next, we prove the following (cf.\ \cite{PanWoo-sd} and the references therein):

\begin{prop} \label{prop:pullbackmono}
Let $(M^4,c_M)$ be a four-dimensional oriented conformal manifold and let $(N^3,c_N,D^N)$
be a three-dimensional Weyl space; denote by $\t_M$ and $\t_N$ the almost twistorial structures
of Examples \ref{exm:smoothtwiststr4} and \ref{exm:smoothtwiststr3}\,, associated to
$(M^4,c_M)$ and $(N^3,c_N,D^N)$\,, respectively.\\
\indent
Let $P$ be a principal bundle over $N^3$ endowed with a connection $\G$ and let $A$ be a nowhere zero
section of $L^*\otimes{\rm Ad}P$\,. Also, let $\phi:(M^4,c_M)\to(N^3,c_N)$ be a surjective horizontally
conformal submersion with connected fibres.\\
\indent
Then any two of the following assertions imply the third:\\
\indent
\quad{\rm (i)} $\phi:(M^4,\t_M)\to(N^3,\t_N)$ is twistorial.\\
\indent
\quad{\rm (ii)} $(A,\G)$ is a monopole on $(N^3,c_N,D^N)$\,.\\
\indent
\quad{\rm (iii)} $\widetilde{\G}$ is anti-self-dual.
\end{prop}
\begin{proof}
Let $\widetilde{R}$ be the curvature form of $\widetilde{\G}$. A straightforward calculation shows
that, up to an anti-self-dual term, the following equality holds
\begin{equation} \label{e:pullbackmono}
\widetilde{R}=\phi^*(R)+\phi^*\bigl((D^N\otimes\nabla\bigr)(A)\bigr)+\bigl(\phi^*(D^N)-\H D_+\bigr)\wedge\phi^*(A)
\end{equation}
where we have used the isomorphism $\Lambda^2(T^*M)=\Lambda^2\bigl(\phi^*(T^*N)\bigr)\oplus\bigl(\phi^*(TN)\otimes\V^*\bigr)$\,,
induced by $c_M$\,. The proof follows.
\end{proof}

\begin{rem}
The fact that \eqref{e:pullbackmono} holds, up to an anti-self-dual term, does not require $\phi$ be horizontally conformal;
moreover, this relation characterizes $\H D_+$ among the partial connections on $\V$, over $\H$.\\
\indent
It follows that \emph{a submersion from a four-dimensional oriented conformal manifold to a three-dimensional Weyl space
is twistorial, as in Example \ref{exm:smoothtwistmap4to3}\,, if and only if it pulls-back any (local) monopole to an anti-self-dual
connection}; the `only if' part is essentially due to \cite{PanWoo-sd} whilst the `if' part is an immediate consequence of
Proposition \ref{prop:pullbackmono} (see also \cite{PanWoo-sd} and the references therein).
\end{rem}

\section{Twistorial maps and harmonic morphisms} \label{section:twistmapsharmorphs}

\indent
We start this section by recalling, from \cite{LouPan}\,, the relations between harmonic morphisms and
the twistorial maps of Examples \ref{exm:smoothtwistmap2to2}\,, \ref{exm:smoothtwistmap3to2}\,,
\ref{exm:smoothtwistmap4to2} and \ref{exm:smoothtwistmap4to3}\,.
Firstly, any map between two-dimensional orientable conformal
manifolds is a harmonic morphism if and only if, with respect to suitable orientations,
it is a twistorial map.

\begin{prop}[\,\cite{LouPan}\,] \label{prop:twistharmorph}
Let $(M^m,c_M,D^M)$ and $(N^n,c_N,D^N)$ be Weyl spaces of dimensions $m$ and $n$\,,
respectively, where $(m,n)\in\{(3,2),(4,2),(4,3)\}$\,. If $m$ $(n)$ is even then
$M^m$ $(N^n)$ is assumed to be oriented. Endow $M^m$ and $N^n$, according to their
dimensions, with the almost twistorial structures $\t_M$ and $\t_N$, respectively,
of Examples \ref{exm:smoothtwiststr2}\,, \ref{exm:smoothtwiststr3} or
\ref{exm:smoothtwiststr4}\,. Let $\phi:M^m\to N^n$ be a submersion.\\
\indent
{\rm (i)} If $(m,n)=(3,2)$ then the following assertions are equivalent:\\
\indent
\quad{\rm (i1)} $\phi:(M^3,c_N,D^N)\to(N^2,c_N,D^N)$ is a harmonic morphism.\\
\indent
\quad{\rm (i2)} $\phi:(M^3,\t_M)\to(N^2,\t_N)$ is a twistorial map.\\
\indent
{\rm (ii)} If $(m,n)=(4,2)$ then the following assertions are equivalent:\\
\indent
\quad{\rm (ii1)} $\phi:(M^4,c_N,D^N)\to(N^2,c_N,D^N)$ is a harmonic morphism
and $\phi:(M^4,\t_M)\to(N^2,\t_N)$ is a twistorial map.\\
\indent
\quad{\rm (ii2)} $\phi:(M^4,c_M)\to(N^2,c_N)$ is horizontally conformal and
the fibres of $\Phi$ are tangent to the connection induced by $D^M$ on $P_M$.\\
\indent
{\rm (iii)} If $(m,n)=(4,3)$ then any two of the following assertions
imply the third:\\
\indent
\quad{\rm (iii1)} $\phi:(M^4,c_N,D^N)\to(N^3,c_N,D^N)$ is a harmonic morphism.\\
\indent
\quad{\rm (iii2)} $\phi:(M^4,\t_M)\to(N^3,\t_N)$ is a twistorial map.\\
\indent
\quad{\rm (iii3)} The fibres of $\Phi$ are tangent to the connection induced by
$D^M$ on $P_M$.
\end{prop}

\indent
Next we discuss twistorial maps between manifolds endowed with the almost twistorial
structures of Examples \ref{exm:smoothtwiststr2}\,, \ref{exm:smoothtwiststr3'} or
\ref{exm:smoothtwiststr4'}\,. The notations are as in Section \ref{twistmaps}\,.\\
\indent
The following proposition is essentially due to \cite{EelSal}\,; its proof is similar
to the proofs of Propositions \ref{prop:smoothtwistmap3'to2} and \ref{prop:smoothtwistmap4'to2}\,.

\begin{prop} \label{prop:smoothtwistmap2to3'4'}
Let $(M^m,c,D)$ be a Weyl space, $m=3,4$\,, and let $k$ be a section of the dual
of the line bundle of $M^m$; if $m=4$ assume $M^4$ oriented and $k=0$\,.
Let $\t_M'$ be the almost twistorial structure on $M^m$ given by Examples \ref{exm:smoothtwiststr3'}
or \ref{exm:smoothtwiststr4'} according to $m=3$ or $m=4$\,, respectively.\\
\indent
Let $N^2$ be an oriented surface in $M^m$. Let $\t_N$ be the twistorial structure
of Example \ref{exm:smoothtwiststr4}\,, associated to $(N^2,c|_N)$\,.\\
\indent
The following assertions are equivalent:\\
\indent
\quad{\rm (i)} $(N^2,\t_N)\to(M^m,\t_M')$ is twistorial (with respect to $p_N$).\\
\indent
\quad{\rm (ii)} $N^2$ is a minimal surface in $(M^m,c,D)$ and $k|_N=0$\,.
\end{prop}

\indent
Let $\phi:(M^3,c)\to N^2$ be a submersion from a three-dimensional conformal manifold
to an oriented two-dimensional manifold. Let $\V={\rm ker \dif\!\phi}$ and $\H=\V^{\perp}$.
{}From the fact that $N^2$ is oriented it follows that $c$ induces an isomorphism
between $\V$ and the line bundle of $M^3$\,.\\
\indent
Let $B^{\H,D}$ be the second fundamental form of $\H$\,, with respect to $D$\,,
and let $I^{\H}$ be the integrability tensor of $\H$\,. We shall denote by $*_{\H}$
the Hodge star operator of $(\H,c_M|_{\H})$\,.

\begin{prop} \label{prop:smoothtwistmap3'to2}
Let $(M^3,c_M,D)$ be a three-dimensional Weyl space and let $k$ be a section of the dual of the line
bundle of $M^3$. Let $\t_M'$ be the almost twistorial structure of Example \ref{exm:smoothtwiststr3'}\,,
associated to $(M^3,c_M,D,k)$\,; denote by $\nabla$ the connection associated to $D$ and $k$\,.\\
\indent
Let $(N^2,c_N)$ be a two-dimensional oriented conformal manifold. Let $\t_N$ be the
twistorial structure of Example \ref{exm:smoothtwiststr2}\,, associated to $(N^2,c_N)$\,.\\
\indent
Let $\phi:M^3\to N^2$ be a submersion. The following assertions are equivalent:\\
\indent
\quad{\rm (i)} $\phi:(M^3,\t_M')\to(N^2,\t_N)$ is a twistorial map.\\
\indent
\quad{\rm (ii)} $\phi:(M^3,c_M)\to(N^2,c_N)$ is horizontally conformal, $D$ is the
Weyl connection of $(M^3,c_M,\V)$ and $*_{\H}I^{\H}=k$\,.\\
\indent
\quad{\rm (iii)} $\phi:(M^3,c_M,D)\to(N^2,c_N)$ is a harmonic morphism, $*_{\H}I^{\H}=k$
and $\trace_{c_M}(B^{\H,D})=0$\,.\\
\indent
\quad{\rm (iv)} $\phi:(M^3,c_M)\to(N^2,c_N)$ is a horizontally conformal submersion and
$\trace_{c_M}(\nabla\!F^{\phi})=0$\,, $\trace_{c_M}(B^{\H,D})=0$\,.
\end{prop}
\begin{proof}
Let $p_{\phi}$ be the section of $P$ corresponding to $F^{\phi}$ where $\t_M'=(P,M,\p,\Fa')$
(see Example \ref{exm:smoothtwistmap3to2}\,). Then condition (1) of
Definition \ref{defn:smoothtwistmap} is satisfied by $\phi$\,, with respect
to $\Phi=\phi\circ{p_{\phi}}^{-1}$\,, if and only if
$\dim_{\C}\bigl(\F'_{p_{\phi}(x)}\cap\dif\!p_{\phi}(T_x^{\C}M)\bigr)=2$\,, for all $x\in M$,
where $\F'$ is the direct sum of the eigenbundles of $\Fa'$ corresponding to $0$ and $-{\rm i}$\,.
Equivalently, condition (1) of Definition \ref{defn:smoothtwistmap} is satisfied by $\phi$\,,
with respect to $\Phi$\,, if and only if $p_{\phi}:(M^3,F^{\phi})\to(P^5,\Fa')$ is holomorphic.\\
\indent
It follows that assertion (i) holds if and only if $p_{\phi}:(M^3,F^{\phi})\to(P^5,\Fa')$
and $\phi:(M^3,F^\phi)\to(N^2,J^N)$ are holomorphic where $J^N$ is the positive
Hermitian structure on $(N^2,c_N)$\,. The latter condition is equivalent to
$\phi:(M^3,c_M)\to(N^2,c_N)$ horizontally conformal.\\
\indent
We can find local representatives $g=U\odot U+2Y\odot\overline{Y}$ of $c_M$ such that
$U$ and $Y$ (locally) generate the eigenbundles of $F^{\phi}$
corresponding to $0$ and $-{\rm i}$\,, respectively; in particular,
$\dif\!\phi(U)=0$\,. Hence, under the identification of $P$ with the bundle of two-dimensional
degenerate spaces on $(M^3,c_M)$\,, the section $p_{\phi}$ corresponds to the
distribution (locally) generated by $U$ and $Y$.\\
\indent
From Proposition \ref{prop:autoparallel} it follows that $p_{\phi}:(M^3,F^{\phi})\to(P^5,\Fa')$
is holomorphic if and only if $g(D_UU,Y)=0$ and $g(D_{\overline{Y}}U,Y)=0$\,. Thus,
after a straightforward calculation we obtain that $p_{\phi}:(M,F^{\phi})\to(P,\Fa'\,)$
is holomorphic if and only if $*_{\H}I^{\H}=k$\,, $\trace_{c_M}(B^{\H,D})=0$ and
the fibres of $\phi$ are geodesics of $D$\,.\\
\indent
Thus, we have proved that (i)$\iff$(ii)$\iff$(iii)\,.\\
\indent
Also, if $\phi:(M^3,c_M)\to(N^2,c_N)$ is horizontally conformal then we have
$\trace_{c_M}(\nabla\!F^{\phi})=0$ if and only if $*_{\H}I^{\H}=k$ and the fibres
of $\phi$ are geodesics of $D$\,. Hence, (iii)$\iff$(iv)\,.
\end{proof}

\indent
Also, we have the following result.

\begin{prop} \label{prop:smoothtwistmap4'to2}
Let $(M^4,c_M,D)$ be a four-dimensional oriented Weyl space. Let $\t_M'$ be the almost twistorial
structure of Example \ref{exm:smoothtwiststr4'}\,, associated to $(M^4,c_M,D)$\,.\\
\indent
Let $(N^2,c_N)$ be a two-dimensional oriented conformal manifold. Let $\t_N$ be the
twistorial structure of Example \ref{exm:smoothtwiststr2}\,, associated to $(N^2,c_N)$\,.\\
\indent
Let $\phi:M^4\to N^2$ be a submersion. The following assertions are equivalent:\\
\indent
\quad{\rm (i)} $\phi:(M^4,\t_M')\to(N^2,\t_N)$ is a twistorial map.\\
\indent
\quad{\rm (ii)} $\phi:(M^4,c_M)\to(N^2,c_N)$ is horizontally conformal and $D$ is the
Weyl connection of $(M^4,c_M,J^{\phi})$.\\
\indent
\quad{\rm (iii)} $\phi:(M^4,c_M,D)\to(N^2,c_N)$ is a harmonic morphism and we have
$\trace_{c_M}(B^{\H,D})=J^{\phi}(*_{\H}I^{\H})$\,.
\end{prop}
\begin{proof}
Let $p_{\phi}$ be the section of $P$ corresponding to $J^{\phi}$ where $\t_M'=(P,M,\p,\J')$
(see Example \ref{exm:smoothtwistmap4to2}\,). Similarly to the proof of
Proposition \ref{prop:smoothtwistmap4'to2}\,, assertion (i) holds if and only if
$p_{\phi}:(M^4,J^{\phi})\to(P^6,\J')$ is holomorphic and
$\phi:(M^4,c_M)\to(N^2,c_N)$ is horizontally conformal. It follows that (i)$\iff$(ii)\,.\\
\indent
We can find local representatives $g=2\bigl(U\odot\overline{U}+Y\odot\overline{Y}\bigr)$
of $c_M$ such that $U$ and $Y$ (locally) generate the eigenbundle of $J^{\phi}$
corresponding $-{\rm i}$ and such that $U$ and $\overline{U}$ are vertical.
Hence, under the identification of $P$ with the bundle of
self-dual spaces on $(M^4,c_M)$\,, the section $p_{\phi}$ corresponds to the
distribution (locally) generated by $U$ and $Y$.\\
\indent
From Proposition \ref{prop:autoparallel} it follows that $p_{\phi}:(M^4,J^{\phi})\to(P^6,\J')$
is holomorphic if and only if $g(D_{\overline{U}}U,Y)=0$ and $g(D_{\overline{Y}}Y,U)=0$\,.
Now, $g(D_{\overline{U}}U,Y)=0$ if and only if the fibres of $\phi$ are minimal with respect
to $D$ whilst $g(D_{\overline{Y}}Y,U)=0$ if and only if
$\trace_{c_M}(B^{\H,D})=J^{\phi}(*_{\H}I^{\H})$\,. Hence, (i)$\iff$(iii)\,.
\end{proof}

\section{Harmonic morphisms and twistorial maps between Weyl spaces of dimensions four and three} \label{section:twistharmorph4to3}

\indent
We start this section by introducing a generalization of the almost twistorial structure of
Example \ref{exm:smoothtwiststr3'}\,.

\begin{exm} \label{exm:smoothtwiststr3''}
Let $(M^3,c)$ be a three-dimensional conformal manifold endowed with two Weyl connections
$D'$ and $D''$. Let $k$ be a section of $L^*$ and let $\nabla$ be the connection
associated to $D''$ and $k$.\\
\indent
Let $\p:P\to M$ be the bundle of nonzero skew-adjoint $f$-structures on $(M^3,c)$\,.
We denote by $\H^0$ the subbundle of $T^{\C}\!P$ such that, at each $p\in P$,
the subspace $\H^0_p\subseteq T^{\C}_p\!P$ is the horizontal lift, with respect to $D'$, of
the eigenspace of $p^{\C}\in{\rm End}(T_{\p(p)}^{\C}\!M)$ corresponding to the eigenvalue $0$\,.
Also, we denote by $\H^{1,0}$ the subbundle of $T^{\C}\!P$ such that, at each $p\in P$,
the subspace $\H^{1,0}_p\subseteq T^{\C}_p\!P$ is the horizontal lift, with respect to $\nabla$,
of the eigenspace of $p^{\C}\in{\rm End}(T_{\p(p)}^{\C}\!M)$ corresponding to the eigenvalue
${\rm i}$\,.\\
\indent
We define the almost $f$-structure $\Fa''$ on $P$ with respect to which
$T^0P=\H^0$ and $T^{1,0}P=({\rm ker}\dif\!\p)^{0,1}\oplus\H^{1,0}$.
Then $\t''=(P,M,\p,\Fa'')$ is a nonintegrable almost twistorial structure on $M$
(the nonintegrability of $\t''$ follows easily from the proof of the integrability result
presented in \cite{Pan-sp}\,).
\end{exm}

\indent
The results of Section \ref{section:twistmapsharmorphs}\,, which involve the almost twistorial
structure of Example \ref{exm:smoothtwiststr3'}\,, can be easily generalized by working, instead,
with the almost twistorial structure of Example \ref{exm:smoothtwiststr3''}\,.
In Theorem \ref{thm:smoothtwistmap4'to3''}\,, below, we shall see that this generalization is
significant when working with maps between Weyl spaces of dimensions four and three.\\
\indent
We shall also need the following.

\begin{exm} \label{exm:smoothtwiststr4'phi}
Let $(M^4,c,D)$ be a four-dimensional oriented Weyl space. Let $\p:P\to M$ be the
bundle of positive orthogonal complex structures on $(M^4,c)$\,.\\
\indent
Let $\phi:M^4\to N^3$ be a submersion; denote by $\V={\rm ker}\dif\!\phi$ and $\H=\V^{\perp}$.
For any $p\in P$ let $l'_p$ and $l''_p$ be the lines in $T^{\C}_{\p(p)}\!M$ spanned by $V-{\rm i}\,p(V)$
and $X-{\rm i}\,p(X)$\,, respectively, for any $V\in\V_{\p(p)}$ and $X\in p(\V_{\p(p)})^{\perp}\cap\H_{\p(p)}$.\\
\indent
We define the almost CR-structure $\mathscr{C}_{\phi}'$ on $P$ which, at each $p\in P$, is the direct sum
of $({\rm ker}\dif\!\p_p)^{0,1}$ and the horizontal lift, with respect to $D$\,, of $l'_p$\,.\\
\indent
Similarly, we define the almost CR-structure $\mathscr{C}_{\phi}''$ on $P$ which, at each $p\in P$,
is the direct sum of $({\rm ker}\dif\!\p_p)^{0,1}$ and the horizontal lift, with respect to $D$\,, of $l''_p$\,.\\
\indent
Then $\t_{\phi}'=(P,M,\p,\mathscr{C}_{\phi}')$ and $\t_{\phi}''=(P,M,\p,\mathscr{C}_{\phi}'')$ are nonintegrable
almost twistorial structures on $M^4$ (the nonintegrability of $\t_{\phi}'$ and $\t_{\phi}''$ follows easily
from the proof of the integrability result presented in \cite{Pan-sp}\,).
\end{exm}

\begin{rem} \label{rem:horconf4to3}
Let $(N^3,c_N)$ be a three-dimensional conformal manifold. At least locally, we may assume that
$(L^*\otimes TN,c)$ is the adjoint bundle of a rank two complex vector bundle $E$ with group ${\rm SU}(2)$\,.
Furthermore, under this identification, $P_N$ is isomorphic to $PE$\,.
As ${\rm SO}(3,\R)={\rm SU}(2)/\mathbb{Z}_2={\rm PSU}(2)$\,, any connection on $(L^*\otimes TN,c)$
corresponds to a connection on the ${\rm PSU}(2)$-bundle $P_N$. Also, $(L^*\otimes TN,c)$
is the adjoint bundle of the ${\rm PSU}(2)$-bundle $P_N$.\\
\indent
Similarly, if $(M^4,c_M)$ is a four-dimensional oriented conformal manifold then the bundle $P_M$ of
positive orthogonal complex structures on $(M^4,c_M)$ is a ${\rm PSU}(2)$-bundle.\\
\indent
Furthermore, a submersion $\phi:(M^4,c_M)\to(N^3,c_N)$ is horizontally conformal
if and only if $\Phi:P_M\to P_N$ is a morphism of ${\rm PSU}(2)$-bundles (in general,
$\Phi$ is a morphism of ${\rm SL}(3,\R)$-bundles). Also, note that, if $\phi$
is horizontally conformal then $P_M=\phi^*(P_N)$ as ${\rm PSU}(2)$-bundles.
\end{rem}

\indent
The following result shows the importance of the almost twistorial structures of
Example \ref{exm:smoothtwiststr4'phi}\,.

\begin{thm} \label{thm:smoothtwistmap4'phi3'}
Let $(M^4,c_M,D^M)$ be a four-dimensional oriented Weyl space and let $P_M$ be the bundle of positive orthogonal
complex structures on $(M^4,c_M)$\,.\\
\indent
Let $(N^3,c_N,D^N)$ be a three-dimensional Weyl space
and let $k$ be a section of the dual of the line bundle of $N^3$; denote by $\nabla$ the connection associated
to $D^N$ and $k$\,. Let $\t_N'$ be the almost twistorial structure
of Example \ref{exm:smoothtwiststr3'}\,, associated to $(N^3,c_N,D^N,k)$\,.\\
\indent
Let $\phi:M^4\to N^3$ be a submersion; denote by $\V={\rm ker}\dif\!\phi$ and $\H=\V^{\perp}$. Let
$\t_{\phi}'$ and $\t_{\phi}''$ be the almost twistorial structures of Example \ref{exm:smoothtwiststr4'phi}\,,
associated to $(M^4,c_M,D^M)$ and $\phi$\,.\\
\indent
{\rm (i)} The following assertions are equivalent:\\
\indent
\quad{\rm (i1)} $\phi:(M^4,\t_{\phi}')\to(N^3,\t_N')$ is twistorial.\\
\indent
\quad{\rm (i2)} $\phi:(M^4,c_M,D^M)\to(N^3,c_N,D^N)$ is a harmonic morphism.\\
\indent
{\rm (ii)} The following assertions are equivalent:\\
\indent
\quad{\rm (ii1)} $\phi:(M^4,\t_{\phi}'')\to(N^3,\t_N')$ is twistorial.\\
\indent
\quad{\rm (ii2)} $\phi:(M^4,c_M)\to(N^3,c_N)$ is horizontally conformal, $\V(D^M-D)=\tfrac12\,k$ and
$\phi^*(D^N)=\H D^M+\tfrac12*_{\H}\!I^{\H}$ as partial connections, over $\H$,
where $D$ is the Weyl connection of $(M^4,c_M,\V)$\,.\\
\indent
\quad{\rm (ii3)} $\phi:(M^4,c_M)\to(N^3,c_N)$ is horizontally conformal
and the partial connections on $P_M$, over $\H$, induced by $D^M$ and $\phi^*(\nabla)$ are equal.
\end{thm}
\begin{proof}
(i) By Remark \ref{rem:horconf4to3}\,, if (i1) holds then $\phi:(M^4,c_M)\to(N^3,c_N)$ is horizontally conformal.\\
\indent
Let $\t_{\phi}'=(P_M,M,\p_M,\mathscr{C}_{\phi}')$ and let $\t'_N=(P_N,N,\p_N,\Fa')$\,.
By associating to each orthogonal complex structure on $(M^4,c_M)$ its eigenspace corresponding to $-{\rm i}$\,,
we identify $\p_M:P_M\to M$ with the bundle of self-dual spaces on $(M^4,c_M)$\,.\\
\indent
Similarly, by associating to each skew-adjoint $f$-structure on $(N^3,c_N)$ the sum of its eigenspaces corresponding
to $0$ and $-{\rm i}$\,, we identify $\p_N:P_N\to N$ with the bundle of (complex) two-dimensional degenerate
spaces on $(N^3,c_N)$\,.\\
\indent
Under this identifications the natural lift $\Phi:P_M\to P_N$ of $\phi$ is given by $\Phi(p)=\dif\!\phi(p)$\,,
for any $p\in P_M$\,.\\
\indent
Let $q$ be a local section of $P_N$\,, over some open set $V\subseteq N$, and let $p$ be the local section of $P_M$\,,
over $\phi^{-1}(V)$\,, such that $\Phi\circ p=q\circ\phi$\,. At least locally, we may assume that $p$ is generated by $Y=U+{\rm i}X$
and $Z$, with $U$ vertical, and $X$ and $Z$ basic. Thus, $q$ is generated by $\dif\!\phi(X)$ and $\dif\!\phi(Z)$\,.\\
\indent
We shall assume that, for some $x_0\in\phi^{-1}(V)$\,, we have $\dif\!q(\dif\!\phi(X_{x_0}))$ horizontal, with respect
to the connection induced by $\nabla$ on $P_N$\,. Then Proposition \ref{prop:autoparallel}\,, the fundamental equation
and a straightforward calculation show that $\dif\!p(\overline{Y}_{x_0})\in{(\mathscr{C}_{\phi}')}_{p(x_0)}$
if and only if $c_M\bigl(\trace_{c_M}(\widetilde{D}\!\dif\!\phi),Z\bigr)=0$\,,
where $\widetilde{D}$ is the connection induced by $D^M$ and $D^N$ on $\phi^*(TN)\otimes T^*M$.\\
\indent
This proves that (i2)$\Longrightarrow$(i1)\,.\\
\indent
To prove the converse, assume that (i1) holds. Then there exists a unique
$A\in{(\mathscr{C}_{\phi}')}_{p(x_0)}$ such that $\dif\!\Phi(A)=\dif\!q(\dif\!\phi(X_{x_0}))$\,.
From the fact that $\dif\!\p_M\bigl({(\Cal_{\phi}')}_{p(x_0)}\bigr)$ is the line spanned by $\overline{Y}_{x_0}$\,,
it follows quickly that $\dif\!\p_M(A)=-\overline{Y}_{x_0}$\,. Together with
$\dif\!\Phi(A)=\dif\!\Phi(\dif\!p(-\overline{Y}_{x_0}))$ this gives $A=\dif\!p(-\overline{Y}_{x_0})$\,.
Hence, $\dif\!p(\overline{Y}_{x_0})\in{(\mathscr{C}_{\phi}')}_{p(x_0)}$ and, therefore,
$c_M\bigl(\trace_{c_M}(\widetilde{D}\!\dif\!\phi),Z\bigr)=0$\,.\\
\indent
(ii) By Remark \ref{rem:horconf4to3}\,, if (ii1) holds then $\phi:(M^4,c_M)\to(N^3,c_N)$ is
horizontally conformal.\\
\indent
Similarly to above, let $p$ be a basic local section of $P_M$\,; that is, there exists
a local section $q$ of $P_N$ such that $\Phi\circ p=q\circ\phi$\,.
At least locally, we may assume that $p$ is generated by $Y=U+{\rm i}X$
and $Z$, with $U$ vertical, and $X$ and $Z$ basic. Thus, $q$ is generated by
$\dif\!\phi(X)$ and $\dif\!\phi(Z)$\,.\\
\indent
To prove (ii1)$\iff$(ii2)\,, we shall asume that $p$ is horizontal,
at some point $x_0\in M$, with respect
to the connection induced by $D^M$ on $P_M$ (this could be done
as follows: firstly, define $p$ over some hypersurface, containing $x_0$ and
which is transversal to the fibres of $\phi$\,, such that $p$ to be horizontal at $x_0$\,;
then, extend $p$\,, to an open neighbourhood of $x_0$\,, so that to be basic); in particular,
$\dif\!p(\overline{Z}_{x_0})\in{(\mathscr{C}_{\phi}'')}_{p(x_0)}$\,.\\
\indent
Assertion (ii1) is equivalent to the fact that, for any such $p$\,, we have
$\dif\!q(\dif\!\phi(X_{x_0}))$ is contained in the eigenbundle of $\Fa'$
corresponding to the eigenvalue ${\rm i}$\,. By using
Proposition \ref{prop:autoparallel}\,, it quickly follows that (ii1)$\iff$(ii2)\,.\\
\indent
Let $h$ be an oriented representative of $c_N$ and let $g$ be a representative of $c_M$
such that $\phi:(M^4,g)\to(N^3,h)$ is a Riemannian submersion. Let $\a^M$ and $\a^N$
be the Lee forms of $D^M$ and $D^N$ with respect to $g$ and $h$\,,
respectively.\\
\indent
The equivalence (ii2)$\iff$(ii3) follows from the fact that any two of the following
assertions imply the third:\\
\indent
\quad(a) $\dif\!p(\H_{x_0})$ is horizontal, with respect to the connection
induced by $D^M$ on $P_M$\,.\\
\indent
\quad(b) $q$ is horizontal at $\phi(x_0)$\,, with respect to the connection
induced by $\nabla$ on $P_N$\,.\\
\indent
\quad(c) At $x_0$ we have $\a^M|_{\V}=\tfrac12\,k$ and
$\a^N=\a^M|_{\H}+\tfrac12*_{\H}\!I^{\H}$\,, where we have identified $k$ and $\a^N$
with their pull-backs by $\phi$ and, in the first equality, we have used the
isomorphism between $\V$ and the pull-back by $\phi$ of the line bundle of $N^3$
induced by $\phi$ and the orientation of $M^4$.\\
\indent
The proof is complete.
\end{proof}

\begin{rem}
Let $(M^4,c,D)$ be a four-dimensional oriented Weyl spaces and let $P$ be the bundle
of positive orthogonal complex structures on $(M^4,c)$\,.\\
\indent
With the same notations as in Example \ref{exm:smoothtwiststr4'phi}\,, by
Theorem \ref{thm:smoothtwistmap4'phi3'}(i)\,, \emph{the relation
$\phi\longleftrightarrow\mathscr{C}_{\phi}'$ induces a bijective correspondence
between one-dimensional foliations on $(M^4,c,D)$\,, which are locally defined
by harmonic morphisms, and certain almost CR-structures on $P$.}
\end{rem}

\indent
Next, we prove the following.

\begin{thm} \label{thm:smoothtwistmap4'to3''}
Let $(M^4,c_M,D^M)$ be a four-dimensional oriented Weyl space and
let $\t_M'$ be the almost twistorial structure
of Example \ref{exm:smoothtwiststr4'}\,, associated to $(M^4,c_M,D^M)$\,.\\
\indent
Let $(N^3,c_N)$ be a three-dimensional conformal manifold endowed with two Weyl connections
$D'$ and $D''$. Let $k$ be a section of the dual of the line bundle of $N^3$ and let $\t_N''$
be the almost twistorial structure of Example \ref{exm:smoothtwiststr3''}\,,
associated to $(N^3,c_N,D',D'',k)$\,.\\
\indent
Let $\phi:M^4\to N^3$ be a submersion. The following assertions are equivalent:\\
\indent
\quad{\rm (i)} $\phi:(M^4,\t_M')\to(N^3,\t_N'')$ is twistorial.\\
\indent
\quad{\rm (ii)} $\phi:(M^4,c_M)\to(N^3,c_N)$ is horizontally conformal and the connection induced by
$D^M$ on the bundle of positive orthogonal complex structures on $(M^4,c_M)$ is the pull-back by
$\phi$ of $(A,\nabla)$\,, where $A=(D''-D')^{\sharp_{c_N}}$ and $\nabla$ is the connection associated
to $D''$ and $k$\,.\\
\indent
\quad{\rm (iii)} The following assertions hold:\\
\indent
\quad\quad{\rm (iii1)} $\phi:(M^4,c_M,D^M)\to(N^3,c_N,D')$ is a harmonic morphism;\\
\indent
\quad\quad{\rm (iii2)} $\phi:(M^4,\t_M)\to(N^3,\t_N)$ is twistorial, where $\t_M$ and $\t_N$ are the almost twistorial structures
of Examples \ref{exm:smoothtwiststr4} and \ref{exm:smoothtwiststr3}\,, associated to $(M^4,c_M)$
and $(N^3,c_N,2D''-D')$\,,
respectively.\\
\indent
\quad\quad{\rm (iii3)} $\V(D^M-D)=\tfrac12\,k$ where $\V={\rm ker}\dif\!\phi$ and
$D$ is the Weyl connection of $(M^4,c_M,\V)$\,.
\end{thm}
\begin{proof}
From Theorem \ref{thm:smoothtwistmap4'phi3'}\,, it follows that assertion (i) holds if and only if
$\phi:(M^4,c_M,D^M)\to(N^3,c_N,D')$ is a harmonic morphism, $\V(D^M-D)=\tfrac12\,k$ and
$\phi^*(D^N)=\H D^M+\tfrac12*_{\H}\!I^{\H}$ as partial connections, over $\H$\,.
Together with the fundamental equation, this quickly gives (i)$\iff$(iii)\,.\\
\indent
Let $\t'_M=(P_M,M,\p_M,\J')$\,. {}From Theorem \ref{thm:smoothtwistmap4'phi3'}\,, it follows that assertion
(i) holds if and only if $\phi:(M^4,c_M)\to(N^3,c_N)$ is horizontally conformal, the partial connections on $P_M$, over $\H$,
induced by $D^M$ and $\phi^*(\nabla)$ are equal, and we have $\H(D^M-D)=D'-D''+\tfrac12*_{\H}\!I^{\H}$,
as partial connections, over $\H$.\\
\indent
Let $(X_1,\ldots,X_4)$ be a positive conformal local frame on $(M^4,c_M)$ such that $X_1$ is vertical and
$X_2$\,, $X_3$\,, $X_4$ are basic. Let $g$ be the local representative of $c_M$ induced by $(X_1,\ldots,X_4)$\,.
Let $\a^M$ and $\a$ be the Lee forms, with respect to $g$\,, of $D^M$ and $D$\,, respectively.
Denote by $\G^i_{jk}$\,, $(i,j,k=1,\ldots,4)$\,, the Christoffel symbols of $D^M$ with respect to $(X_1,\ldots,X_4)$\,.
Then a straightforward calculation gives the following relations:
\begin{equation*}
\begin{split}
\G^1_{21}+\G^3_{41}&=\bigl(\a^M-\a-\tfrac12*_{\H}\!I^{\H}\bigr)(X_2)\;,\\
\G^1_{31}-\G^2_{41}&=\bigl(\a^M-\a-\tfrac12*_{\H}\!I^{\H}\bigr)(X_3)\;,\\
\G^1_{41}+\G^2_{31}&=\bigl(\a^M-\a-\tfrac12*_{\H}\!I^{\H}\bigr)(X_4)\;.
\end{split}
\end{equation*}
\indent
Thus, we have proved that (i) holds if and only if $\phi:(M^4,c_M)\to(N^3,c_N)$ is horizontally conformal,
the partial connections on $P_M$, over $\H$, induced by $D^M$ and $\phi^*(\nabla)$ are equal, and the
following relations hold:
\begin{equation*}
\begin{split}
\G^1_{21}+\G^3_{41}&=(D'-D'')(X_2)\;,\\
\G^1_{31}-\G^2_{41}&=(D'-D'')(X_3)\;,\\
\G^1_{41}+\G^2_{31}&=(D'-D'')(X_4)\;.
\end{split}
\end{equation*}
It follows that (i)$\iff$(ii)\,.
\end{proof}

\begin{rem}
1) Let $(M^4,c_M,D^M)$ be a four-dimensional oriented Weyl space.
Let $\t_M'$ be the almost twistorial structure
of Example \ref{exm:smoothtwiststr4'}\,, associated to $(M^4,c_M,D^M)$\,;
denote by $P_M$ be the bundle of positive orthogonal complex structures on $(M^4,c_M)$\,.\\
\indent
Let $\phi:M^4\to N^3$ be a submersion onto a three-dimensional manifold. Let $P_N$ be the bundle of
oriented lines on $N^3$. Then, similar to Example \ref{exm:smoothtwistmap4to3}\,, there can be
defined a bundle map $\Phi:P_M\to P_N$.\\
\indent
Suppose that there exists an almost twistorial structure $\t$ on $N^3$ such that
$\phi:(M^4,\t_M')\to(N^3,\t)$ is a twistorial map, with respect to $\Phi$.
Then there exists a section $k$ of the dual of the line bundle of $N^3$,
a conformal structure $c_N$ on $N^3$ and Weyl connections $D'$ and $D''$
on $(N^3,c_N)$ such that $\t$ is the almost twistorial structure of
Example \ref{exm:smoothtwiststr3''}\,, associated to $(N^3,c_N,D',D'',k)$\,.\\
\indent
Similar comments apply to all of the previous examples of submersive twistorial maps.\\
\indent
2) Let $(M^4,c_M)$ be an oriented four-dimensional conformal manifold and let $(N^3,c_N,D)$ be a
three-dimensional Weyl space. Denote by $\t_M$ and $\t_N$ the almost twistorial structures,
of Examples \ref{exm:smoothtwiststr4} and \ref{exm:smoothtwiststr3}\,, associated to $(M^4,c_M)$
and $(N^3,c_N,D)$\,, respectively. Let $\phi:(M^4,\t_M)\to(N^3,\t_N)$ be a twistorial map
(see \cite{PanWoo-sd} and \cite{Cal-sds} for examples of such maps).\\
\indent
Let $D'$ be a Weyl connection on $(N^3,c_N)$ and let $k$ be a section of the dual of the line bundle
of $N^3$; denote by $D''=\tfrac12(D+D')$\,.\\
\indent
From Theorem \ref{thm:smoothtwistmap4'to3''}\,, it follows that there exists a unique Weyl connection
$D^M$ on $(M^4,c_M)$ such that $\phi:(M^4,\t_M')\to(N^3,\t_N'')$ is twistorial, where
$\t_M'$ and $\t_N''$ are the almost twistorial structures, of Examples \ref{exm:smoothtwiststr4'} and
\ref{exm:smoothtwiststr3''}\,, associated to $(M^4,c_M,D^M)$ and $(N^3,c_N,D',D'',k)$\,, respectively.\\
\indent
3) Let $(M^4,c_M)$ be a four-dimensional oriented conformal manifold and let $L$ be the line bundle
of $M^4$. Denote by $P_{\pm}$ the bundles of positive/negative orthogonal complex structures on $(M^4,c_M)$\,.
As $P_{\pm}$ are, locally, the projectivisations of the bundles of positive/negative spinors on
$(L^*\otimes TM,c_M)$\,, any conformal connection on $(M^4,c_M)$ corresponds to a pair
$(\G_+,\G_-)$\,, where $\G_{\pm}$ are connections on $P_{\pm}$\,.\\
\indent
Let $\phi:(M^4,c_M)\to(N^3,c_N)$ be a horizontally conformal submersion onto a three-dimensional
conformal manifold. Endow $(N^3,c_N)$ with two Weyl connections $D'$ and $D''$ and let $k$ be a section
of the dual of the line bundle of $N^3$. Let $\G_+$ be the connection on $P_+$ which is the pull-back by
$\phi$ of $(A,\nabla)$\,, where $A=(D''-D')^{\sharp_{c_N}}$ and $\nabla$ is the connection associated
to $D''$ and $k$\,.\\
\indent
Let $\G_-$ be a connection on $P_-$ and suppose that the connection $D^M$ corresponding
to $(\G_+,\G_-)$ is torsion-free. Then, obviously, $D^M$ is a Weyl connection
on $(M^4,c_M)$\,. Moreover, by Theorem \ref{thm:smoothtwistmap4'to3''}\,, the map $\phi:(M^4,\t_M')\to(N^3,\t_N'')$
is twistorial, where $\t_M'$ and $\t_N''$ are the almost twistorial structures, of
Examples \ref{exm:smoothtwiststr4'} and \ref{exm:smoothtwiststr3''}\,, associated to
$(M^4,c_M,D^M)$ and $(N^3,c_N,D',D'',k)$\,, respectively.\\
\indent
4) With the same notations as in Theorem \ref{thm:smoothtwistmap4'to3''}\,, it can be proved
that $(A,\nabla)$ is a monopole on $(N^3,c_N,2D''-D')$ if and only if $D'=D''$ and $\nabla$ is flat
(cf.\ \cite{GauTod}\,; see \cite{LouPan} for details about the resulting maps).
\end{rem}

\indent
The following result is an immediate consequence of Theorem \ref{thm:smoothtwistmap4'to3''}\,;
note that, the equivalence (ii)$\iff$(iii) appears in \cite{LouPan}\,.

\begin{cor} \label{cor:smoothtwistmap4'to3'}
Let $(M^4,c_M,D^M)$ be a four-dimensional oriented Weyl space and let $\t_M'$ be the almost twistorial structure
of Example \ref{exm:smoothtwiststr4'}\,, associated to $(M^4,c_M,D^M)$\,.\\
\indent
Let $(N^3,c_N,D^N)$ be a three-dimensional Weyl space and let $k$ be a section of the dual of the line bundle of $N^3$.
Let $\t_N'$ be the almost twistorial structure of Example \ref{exm:smoothtwiststr3'}\,, associated to $(N^3,c_N,D^N,k)$\,.\\
\indent
Let $\phi:M^4\to N^3$ be a submersion. The following assertions are equivalent:\\
\indent
\quad{\rm (i)} $\phi:(M^4,\t_M')\to(N^3,\t_N')$ is twistorial.\\
\indent
\quad{\rm (ii)} $\phi:(M^4,c_M)\to(N^3,c_N)$ is horizontally conformal and the connection induced by $D^M$ on the bundle
of positive orthogonal complex structures on $(M^4,c_M)$ is the pull-back by $\phi$ of\/ $\nabla$, where
$\nabla$ is the connection associated to $D^N$ and $k$\,.\\
\indent
\quad{\rm (iii)} The following assertions hold:\\
\indent
\quad\quad{\rm (iii1)} $\phi:(M^4,c_M,D^M)\to(N^3,c_N,D^N)$ is a harmonic morphism;\\
\indent
\quad\quad{\rm (iii2)} $\phi:(M^4,\t_M)\to(N^3,\t_N)$ is twistorial, where $\t_M$ and $\t_N$ are the almost twistorial structure
of Examples \ref{exm:smoothtwiststr4} and \ref{exm:smoothtwiststr3}\,, associated to $(M^4,c_M)$ and $(N^3,c_N,D^N)$\,,
respectively.\\
\indent
\quad\quad{\rm (iii3)} $\V(D^M-D)=\tfrac12\,k$ where $\V={\rm ker}\dif\!\phi$ and $D$ is the Weyl connection of $(M^4,c_M,\V)$\,.
\end{cor}

\begin{exm}
1) Let $\phi:(M^4,g)\to(N^3,h)$ be a harmonic morphism given by the Gibbons-Hawking or the Beltrami fields construction
(see \cite{PanWoo-exm}\,). Then $\phi$ satisfies (ii) and (iii) of Corollary \ref{cor:smoothtwistmap4'to3'}\,,
with $D^M$ and $D^N$ the Levi-Civita connections of $g$ and $h$\,, respectively, and a suitable choice of $k$
(\,\cite{LouPan}\,, \cite{PanWoo-sd}\,). Hence, $\phi$ satisfies assertion (i) of Corollary \ref{cor:smoothtwistmap4'to3'}
(and, also, (i) of Theorem \ref{thm:smoothtwistmap4'to3''}\,).\\
\indent
2) Let $\phi:(M^4,g)\to(N^3,h)$ be a harmonic morphism of Killing type between Riemannian manifolds of dimensions
four and three. Then $\phi$ satisfies (iii1)\,, (iii2) of Theorem \ref{thm:smoothtwistmap4'to3''}\,,
with $D^M$ and $D'$ the Levi-Civita connections of $g$ and $h$\,, respectively, and a suitable
choice of $D''$ \cite{PanWoo-sd}\,; also, $\phi$ satisfies (iii3)\,, with $k=0$\,. Hence, $\phi$ satisfies
assertion (i) of Theorem \ref{thm:smoothtwistmap4'to3''}\,.
\end{exm}

\end{document}